# Supersaturation Problem for Color-Critical Graphs

Oleg Pikhurko[a,1], Zelealem B. Yilma[b]

[a]*Mathematics Institute and DIMAP, University of Warwick, Coventry CV4 7AL, UK*
[b]*Carnegie Mellon University Qatar, Doha, QATAR*

**Abstract**

The *Turán function* $\mathrm{ex}(n, F)$ of a graph $F$ is the maximum number of edges in an $F$-free graph with $n$ vertices. The classical results of Turán and Rademacher from 1941 led to the study of supersaturated graphs where the key question is to determine $h_F(n, q)$, the minimum number of copies of $F$ that a graph with $n$ vertices and $\mathrm{ex}(n, F) + q$ edges can have.

We determine $h_F(n, q)$ asymptotically when $F$ is *color-critical* (that is, $F$ contains an edge whose deletion reduces its chromatic number) and $q = o(n^2)$.

Determining the exact value of $h_F(n, q)$ seems rather difficult. For example, let $c_1$ be the limit superior of $q/n$ for which the extremal structures are obtained by adding some $q$ edges to a maximum $F$-free graph. The problem of determining $c_1$ for cliques was a well-known question of Erdős that was solved only decades later by Lovász and Simonovits. Here we prove that $c_1 > 0$ for every color-critical $F$. Our approach also allows us to determine $c_1$ for a number of graphs, including odd cycles, cliques with one edge removed, and complete bipartite graphs plus an edge.

*Keywords:* Extremal Graph Theory, Removal Lemma, Supersaturation, Turán Function

## 1. Introduction

The *Turán function* $\mathrm{ex}(n, F)$ of a graph $F$ is the maximum number of edges in an $F$-free graph with $n$ vertices. In 1907, Mantel [15] proved that $\mathrm{ex}(n, K_3) = \lfloor n^2/4 \rfloor$, where $K_r$ denotes the complete graph on $r$ vertices. The fundamental paper of Turán [24] solved this extremal problem for cliques: the *Turán graph* $T_r(n)$, the complete $r$-partite graph of order $n$ with parts of size $\lfloor n/r \rfloor$ or $\lceil n/r \rceil$, is the unique maximum $K_{r+1}$-free graph of order $n$. Thus we have $\mathrm{ex}(n, K_{r+1}) = t_r(n)$, where $t_r(n) = |E(T_r(n))|$.

Stated in the contrapositive, this implies that a graph with $t_r(n) + 1$ edges (where, by default, $n$ denotes the number of vertices) contains at least one copy of $K_{r+1}$. Rademacher (1941, unpublished) showed that a graph with $\lfloor n^2/4 \rfloor + 1$ edges contains not just one but at least $\lfloor n/2 \rfloor$ copies of a triangle. This is perhaps the first result in the so-called "theory of supersaturated graphs" that focuses on the function

$$h_F(n, q) = \min \left\{ \#F(H) : |V(H)| = n,\ |E(H)| = \mathrm{ex}(n, F) + q \right\},$$

the minimum number of $F$-subgraphs in a graph $H$ with $n$ vertices and $\mathrm{ex}(n, F) + q$ edges. (We say that $G$ is a *subgraph* of $H$ if $V(G) \subseteq V(H)$ and $E(G) \subseteq E(H)$; we call $G$ an *F-subgraph* if it is isomorphic to $F$.) One possible construction is to add some $q$ edges to a maximum $F$-free graph; let $t_F(n, q)$ be the smallest number of $F$-subgraphs that can be achieved this way. Clearly, $h_F(n, q) \leq t_F(n, q)$.

Erdős [3] extended Rademacher's result by showing that $h_{K_3}(n, q) = t_{K_3}(n, q) = q \lfloor n/2 \rfloor$ for $q \leq 3$. Later, he [4, 5] showed that there exists some small constant $\epsilon_r > 0$ such that $h_{K_r}(n, q) = t_{K_r}(n, q)$ for all $q \leq \epsilon_r n$. Lovász and Simonovits [13, 14] found the best possible value of $\epsilon_r$ as $n \to \infty$, settling a

*Email addresses:* O.Pikhurko@warwick.ac.uk (Oleg Pikhurko ), zyilma@qatar.cmu.edu (Zelealem B. Yilma)
[1]Supported by ERC grant 306493 and EPSRC grant EP/K012045/1.

*Preprint submitted to Elsevier*   November 7, 2016

long-standing conjecture of Erdős [3]. If fact, the second paper [14] completely solved the $h_{K_r}(n,q)$-problem when $q = o(n^2)$. The case $q = \Omega(n^2)$ of the supersaturation problem for cliques has been actively studied and proved notoriously difficult. Only recently was an asymptotic solution found: by Razborov [18] for $K_3$ (see also Fisher [8]), by Nikiforov [17] for $K_4$, and by Reiher [19] for general $K_r$.

If $F$ is bipartite, then there is a beautiful (and still open) conjecture of Erdős–Simonovits [22] and Sidorenko [20] whose positive solution would determine $h_F(n,q)$ asymptotically for $q = \Omega(n^2)$. We refer the reader to some recent papers on the topic, [2, 9, 10, 12, 23], that contain many references.

Obviously, if we do not know $\mathrm{ex}(n, F)$, then it is difficult to say much about the supersaturation problem for small $q$. A large and important class of graphs for which the Turán function is well understood is formed by *color-critical* graphs, that is, graphs whose chromatic number can be decreased by removing an edge:

**Definition 1.1.** *A graph $F$ is $r$-critical if $\chi(F) = r + 1$ but $F$ contains an edge $e$ such that $\chi(F - e) = r$.*

Simonovits [21] proved that for an $r$-critical graph $F$ we have $\mathrm{ex}(n, F) = t_r(n)$ for all large enough $n \geq n_0(F)$; furthermore, $T_r(n)$ is the unique maximum $F$-free graph. From now on, we assume everywhere that $F$ is an $r$-critical graph and $n$ is sufficiently large so that, in particular, the above result from [21] applies.

The supersaturation problem for a color-critical graph that is not a clique was first considered by Erdős [7, Page 296] who stated that the methods of [7] can prove that $h_{C_5}(2m, 1) = 2m(2m-1)(2m-2)$, where $C_k$ denotes the cycle of length $k$. Recently, Mubayi [16] embarked on a systematic study of this problem for color-critical graphs:

**Definition 1.2.** *Fix $r \geq 2$ and let $F$ be an $r$-critical graph. Let $c(n, F)$ be the minimum number of copies of $F$ in the graph obtained from $T_r(n)$ by adding one edge.*

Since we assume that $n$ is large, we have that $c(n, F) = t_F(n, 1)$. Also, it is not hard to show that for $q = o(n^2)$ we have
$$qc(n, F) \leq t_F(n, q) \leq (1 + o(1))qc(n, F). \tag{1}$$

**Theorem 1.3** (Mubayi [16]). *For every $r$-critical graph $F$, there exists a constant $c_0 = c_0(F) > 0$ such that for all sufficiently large $n$ and $1 \leq q < c_0 n$ we have*
$$h_F(n, q) \geq qc(n, F). \tag{2}$$

As is pointed out in [16], the bound in (2) is asymptotically best possible. Also, (2) is sharp for some graphs $F$, including odd cycles and $K_4 - e$, the graph obtained from $K_4$ by deleting an edge.

Our Theorems 3.10–3.11 show that in order to determine $h_F(n, q)$ asymptotically for $q = o(n^2)$, it is enough to consider graphs constructed as follows: $V(H) = X \cup V_1 \cup \ldots \cup V_r$ where $|X| = O(q/n)$ and $V_1 \cup \ldots \cup V_r$ span a Turán graph, except $V_1$ contains some extra edges spread uniformly. Determining the asymptotic behavior of $h_F(n, q)$ then reduces to optimizing a function of $|X|$, the neighborhoods of $x \in X$, and the number of extra edges in $V_1$. We solve this problem when $q/n \to \infty$ in Theorem 3.10.

Let $\mathcal{T}_r(n, q)$ be the set of graphs obtained from the Turán graph $T_r(n)$ by adding $q$ edges:
$$\mathcal{T}_r(n, q) = \big\{ H : |V(H)| = n,\ |E(H)| = t_r(n) + q,\ H \supseteq T_r(n) \big\}.$$

These graphs are natural candidates, when $q$ is small, for membership in
$$\mathcal{H}_F(n, q) = \big\{ H : |V(H)| = n,\ |E(H)| = \mathrm{ex}(n, F) + q,\ \#F(H) = h_F(n, q) \big\},$$
the set of graphs on $n$ vertices and $\mathrm{ex}(n, F) + q$ edges which contain the smallest number of copies of $F$. Of particular interest is identifying a threshold for when graphs in $\mathcal{T}_r(n, q)$ are optimal or asymptotically optimal. In view of (1), the threshold constant for the latter property is formally defined as
$$c_2(F) = \sup \big\{ c : \forall \epsilon > 0\ \exists n_0\ \forall n \geq n_0\ \forall q \leq cn\ \ h_F(n, q) \geq (1 - \epsilon)qc(n, F) \big\}. \tag{3}$$

Our Theorem 3.12 determines this parameter for every color-critical $F$. Its statement requires some technical definitions so we postpone it until Section 3. Informally speaking, Theorem 3.12 states that $c_2$ is the limit inferior of $q/n$ when the following construction starts beating the bound $(1-o(1))qc(n,F)$: add a new vertex $x$ of degree $t_r(n) + q - t_r(n-1)$ to $T_r(n-1)$ so that the number of the created $F$-subgraphs is minimized. For some instances of $F$ and values of $q$, this construction indeed wins. On the other hand, there are also examples of $F$ with $c_2(F) = \infty$; in the latter case we prove the stronger claim that $h_F(n,q) = (1+o(1))qc(n,F)$ for all $q = o(n^2)$ (not just for $q = O(n)$), see Theorem 3.12.

We also focus on the optimality of $\mathcal{T}_r(n,q)$ and our result qualitatively extends Theorem 1.3 as follows:

**Theorem 1.4.** *For every $r$-critical graph $F$, there exist $c_1 > 0$ and $n_0$ such that for all $n > n_0$ and $q < c_1 n$, we have $h_F(n,q) = t_F(n,q)$ (in fact, more strongly, we have $\mathcal{H}_F(n,q) \subseteq \mathcal{T}_r(n,q)$).*

A natural question arises here, namely, how large $c_1 = c_1(F)$ in Theorem 1.4 can be. So we define

$$\begin{aligned} q_F(n) &= \max\{q \in \mathbb{N} : h_F(n,q') = t_F(n,q') \text{ for all } q' \leq q\}, \\ c_1(F) &= \liminf_{n \to \infty} \left\{\frac{q_F(n)}{n}\right\}. \end{aligned} \qquad (4)$$

In 1955 Erdős [3] conjectured that $q_{K_3}(n) \geq \lfloor n/2 \rfloor - 1$ and observed that, if true, this inequality would be sharp for even $n$. This conjecture (and even its weaker version if $c_1(K_3) \geq 1/2$) remained open for decades until it was finally proved by Lovász and Simonovits [13, 14] whose more general results imply that $c_1(K_{r+1}) = 1/r$ for every $r$.

Our approach allows us to determine the value of $c_1(F)$ for a number of other graphs. Here are some examples.

**Theorem 1.5.** *Let $F$ be an odd cycle. Then $c_1(F) = 1/2$.*

**Theorem 1.6.** *Let $r \geq 2$ and $F = K_{r+2} - e$ be obtained from $K_{r+2}$ by removing an edge. Then $c_1(F) = (r-1)/r^2$.*

Also, we can determine $c_1(F)$ if $F$ is obtained from a complete bipartite graph by adding an edge (see Corollary 4.8 and Theorem 4.9) and for a whole class of what we call *pair-free* graphs. Unfortunately, these results are rather technical to state so instead we refer the reader to Section 4.

In all these examples (as well as for $F = K_{r+1}$), if $q < (c_1(F) - \epsilon)n$ and $n \geq n_0(\epsilon, F)$, then not only $h_F(n,q) = t_F(n,q)$ but also $\mathcal{H}_F(n,q) \subseteq \mathcal{T}_r(n,q)$, that is, every extremal graph is obtained by adding edges to the Turán graph.

In Theorems 1.5 and 1.6, the $c_1$-threshold coincides with the moment when the number of copies of $F$ may be strictly decreased by using a non-equitable partition. For example, if $F = C_3 = K_3$ and $n = 2\ell$ is even, then instead of adding $q = \ell$ edges to the Turán graph $T_2(n) = K_{\ell,\ell}$, one can add $q+1$ edges to the larger part of $K_{\ell+1,\ell-1}$ and get fewer triangles.

Interestingly, the congruence class of $n$ modulo $r$ may also affect the value of $q_F(n)$, which happens already for $F = K_3$. Indeed, if $n = 2\ell + 1$ is odd and we start with $K_{\ell+2,\ell-1}$ instead of $T_2(n) = K_{\ell+1,\ell}$, then we need to add extra $q + 2$ edges (not $q + 1$ as it is for even $n$); in fact, $q_{K_3}(2\ell+1)$ is about twice as large as $q_{K_3}(2\ell)$. Hence, we also define the following $r$ constants

$$c_{1,i}(F) = \liminf_{\substack{n \to \infty \\ n \equiv i \pmod{r}}} \left\{\frac{q_F(n)}{n}\right\}, \quad 0 \leq i \leq r-1. \qquad (5)$$

Clearly, we have $c_1(F) = \min\{c_{1,i}(F) : 0 \leq i \leq r-1\}$. In some cases, we are able to determine the constants $c_{1,i}(F)$ as well.

Still, some other and more complicated phenomena can occur at the $c_1$-threshold, as the 2-critical graph $F$ with 7 vertices from Section 4.2 demonstrates. Namely, if $n = 2\ell$ is even and we wish to add a small number of extra edges to $T_2(n) = K_{\ell,\ell}$ optimally, then we should divide them equally between the two parts: our $r$-critical graph $F$ was devised so that if its copy uses exactly two of the extra edges, then they belong to the same part. On the other hand, if we start with $K_{\ell+1,\ell-1}$, then it is advantageous to put all extra edges into the larger part, in spite of the greater number of $F$-subgraphs created by such unbalancing. Although these phenomena contribute $O(n^5)$ to $\#F(H)$ (while $h_F(n,q) = \Theta(n^6)$) when $q = \Theta(n)$, the corresponding lower-order terms do affect the value of $c_1(F)$. Also, the proof of Theorem 4.9 shows that some other interesting phenomena occur when $F$ is obtained from $K_{7,2}$ by adding an edge into the larger part and $n$ is odd; however, here $c_1(F)$ is determined by $c_{1,0}(F)$ which behaves 'regularly'. This indicates that a general formula for $c_1(F)$ may be difficult to obtain.

The rest of the paper is organized as follows. In the next section we introduce the functions and parameters with which we work. Our asymptotic results on the case $q = o(n^2)$, including some general lower bounds on $c_1(F)$ as well as the determination of $c_2(F)$, are proved in Section 3. We apply our method to determine $c_1(F)$ for some special graphs in Section 4. And lastly, we have appended a glossary of terms for ease of lookup of the many defined quantities and parameters.

## 2. Parameters

In the arguments and definitions to follow, $F$ will be an $r$-critical graph and we let $f = |V(F)|$ be the number of vertices of $F$. We identify graphs with their edge sets, e.g. $|F|$ means $|E(F)|$. Typically, the order of a graph $H$ under consideration will be denoted by $n$ and viewed as tending to infinity. We will use the asymptotic terminology (such as, for example, the expression $O(1)$) to hide constants independent of $n$. We write $x = y \pm z$ to mean $|x - y| \leq z$. Also, we may ignore rounding errors, when these are not important.

Let us begin with an estimate for $c(n, F)$.

**Lemma 2.1.** *Let $F$ be an $r$-critical graph on $f$ vertices. There is a positive constant $\alpha_F$ such that*

$$c(n, F) = \alpha_F n^{f-2} + O(n^{f-3}).$$

This was proved by Mubayi [16] by providing an explicit formula for $c(n, F)$, see Identity (6) here which we are about to derive.

If $F$ is an $r$-critical graph, we call an edge $e$ (resp., a vertex $v$) a *critical edge* (resp., a *critical vertex*) if $\chi(F - e) = r$ (resp., $\chi(F - v) = r$). Given disjoint sets $V_1, \ldots, V_r$, let $K(V_1, \ldots, V_r)$ be formed by connecting all vertices $v_i \in V_i$ and $v_j \in V_j$ with $i \neq j$, i.e., $K(V_1, \ldots, V_r)$ is the complete $r$-partite graph on vertex classes $V_1, \ldots, V_r$. Let $H$ be obtained from $K(V_1, \ldots, V_r)$ by adding one edge $xy$ into the first part and let $c(n_1, \ldots, n_r; F)$ where $n_i = |V_i|$, denote the number of copies of $F$ contained in $H$. Let $uv \in F$ be a critical edge and let $\chi_{uv}$ be a proper $r$-coloring of $F - uv$ where $\chi_{uv}(u) = \chi_{uv}(v) = 1$. Let $x^i_{uv}$ be the number of vertices of $F$ excluding $u, v$ that receive color $i$. An edge preserving injection of $F$ into $H$ is obtained by picking a critical edge $uv$ of $F$, mapping it to $xy$, then mapping the remaining vertices of $F$ to $H$ so that no two adjacent vertices get mapped to the same part of $H$. Such a mapping corresponds to some coloring $\chi_{uv}$. So, with $\mathrm{Aut}(F)$ denoting the number of automorphisms of $F$, we obtain

$$c(n_1, \ldots, n_r; F) = \frac{1}{\mathrm{Aut}(F)} \sum_{uv \text{ critical}} \sum_{\chi_{uv}} 2\,(n_1 - 2)_{x^1_{uv}} \prod_{i=2}^{r} (n_i)_{x^i_{uv}}, \qquad (6)$$

where $(n)_k = n(n-1)\cdots(n-k+1)$ denotes the *falling factorial*. Lemma 2.1 follows now because $c(n, F)$ is given by (6) for some numbers $n_i = \frac{n}{r} \pm 1$. If $n_1 = \cdots = n_r = n/r$, then (6) is a polynomial in $n$ of degree $f - 2$ (and $\alpha_F$ is the leading coefficient). Also, if $n_1 \leq n_2 \leq \cdots \leq n_r$ and $n_r - n_1 \leq 1$, then $n_i$'s assume at most two different values and we have

$$c(n, F) = \min\{c(n_1, \ldots, n_r; F),\, c(n_r, \ldots, n_1; F)\}. \qquad (7)$$

A recurring argument in our proofs involves moving vertices or edges from one class to another, potentially changing the partition of $n$. To this end, we compare different values of $c(n_1, \ldots, n_r; F)$. In [16], Mubayi proved that $c(n_1, \ldots, n_r; F) \geq c(n, F) + O(an^{f-3})$ for all partitions $n_1 + \ldots + n_r = n$ where $\lfloor n/r \rfloor - a \leq n_i \leq \lceil n/r \rceil + a$ for every $i \in [r]$. We need the following, more precise estimate:

**Lemma 2.2.** *There exist constants $\zeta_F$ and $C_F$ such that the following holds for all large $n$. Let $c(n, F) = c(n'_1, \ldots, n'_r; F)$ as in (7). Let $n_1 + \ldots + n_r = n$. Define $a_i = n_i - n'_i$ for $i \in [r]$ and $A = \max\{|a_i| : i \in [r]\}$. Then*

$$\left|c(n_1, \ldots, n_r; F) - c(n, F) - \zeta_F a_1 n^{f-3}\right| \leq C_F A^2 n^{f-4}.$$

*Proof.* Assume that $A \neq 0$ for otherwise there is nothing to prove. We estimate the value of the polynomial $c(n_1, \ldots, n_r; F)$ using the Taylor expansion about $(n'_1, \ldots, n'_r)$. Namely,

$$c(n'_1 + a_1, \ldots, n'_r + a_r; F) - c(n'_1, \ldots, n'_r; F) - \sum_{j=1}^{r} a_j \frac{\partial c}{\partial_j}(n'_1, \ldots, n'_r), \tag{8}$$

is a polynomial of degree at most $f - 2$ in the variables $n'_i$ and $a_i$ in which every monomial contains at least two $a_i$'s; thus it is $O(A^2 n^{f-4})$. Furthermore, as $|n'_i - n/r| \leq 1$ for all $1 \leq i \leq r$, we have

$$\left|\frac{\partial c}{\partial_i}(n'_1, \ldots, n'_r) - \frac{\partial c}{\partial_i}(n/r, \ldots, n/r)\right| = O(n^{f-4}).$$

Thus, the expression in (8) remains within $O(A^2 n^{f-4})$ if we replace the last sum in (8) by

$$\sum_{j=1}^{r} a_j \frac{\partial c}{\partial_j}(n/r, \ldots, n/r) = a_1 \left(\frac{\partial c}{\partial_1}(n/r, \ldots, n/r) - \frac{\partial c}{\partial_2}(n/r, \ldots, n/r)\right),$$

where we used the facts that, by symmetry, all partial derivatives $\frac{\partial c}{\partial_j}(n/r, \ldots, n/r)$ for $j = 2, \ldots, r$ are equal to each other and that $a_2 + \cdots + a_r = -a_1$. Now, we can let $\zeta_F$ be the coefficient of $n^{f-3}$ in $\frac{\partial c}{\partial_1}(n/r, \ldots, n/r) - \frac{\partial c}{\partial_2}(n/r, \ldots, n/r)$. □

**Definition 2.3.** *For an $r$-critical graph $F$, let $\pi_F = \begin{cases} \frac{\alpha_F}{|\zeta_F|}, & \text{if } \zeta_F \neq 0, \\ \infty, & \text{if } \zeta_F = 0. \end{cases}$*

To give a brief foretaste of the arguments to come, let us compare the number of copies of a 2-critical graph $F$ in some $H \in \mathcal{T}_2(n, q)$ and a graph $H'$ with $K(V_1, V_2) \subseteq H'$ where $n = 2\ell$ is even, $|V_1| = \ell + 1$, and $|V_2| = \ell - 1$. While $H$ contains $q$ 'extra' edges, the identity $(\ell + 1)(\ell - 1) = \ell^2 - 1$ implies that the number of 'extra' edges in $H'$ is $q + 1$. Ignoring, for now, the copies of $F$ that use more than one 'extra' edge, we have to compare the quantities $\#F(H) \approx qc(n, F) \approx q\alpha_F n^{f-2}$ and $\#F(H') \approx (q+1)(\alpha_F n^{f-2} \pm \zeta_F n^{f-3})$. (Note that we can control the sign in front of $\zeta_F$ by choosing the part $V_i$ into which we add all 'extra' edges.) It becomes clear that the ratio $\alpha_F/|\zeta_F|$ will play a significant role in bounding $c_1(F)$.

Another phenomenon of interest is the existence of a vertex with large degree in each part. Let $\boldsymbol{d} = (d_1, \ldots, d_r)$ and let $\#F(n_1, \ldots, n_r; \boldsymbol{d})$ be the number of copies of $F$ in the graph $H = K(V_1, \ldots, V_r) + z$ where $|V_i| = n_i$ and the extra vertex $z$ has $d_i$ neighbors in $V_i$. Let $\#F(n, \boldsymbol{d})$ correspond to the case when $n_1 + \ldots + n_r = n - 1$ satisfy $n_1 \geq \ldots \geq n_r \geq n_1 - 1$.

We have the following formula for $\#F(n_1, \ldots, n_r; \boldsymbol{d})$. An edge preserving injection from $F$ to $H$ is obtained by choosing a critical vertex $u$, mapping it to $z$, then mapping the remaining vertices of $F$ to $H$ so that neighbors of $u$ get mapped to neighbors of $z$ and no two adjacent vertices get mapped to the same part. Such a mapping is given by an $r$-coloring $\chi_u$ of $F - u$. Thus

$$\#F(n_1, \ldots, n_r; \boldsymbol{d}) = \frac{1}{\text{Aut}(F)} \sum_{u \text{ critical}} \sum_{\chi_u} \prod_{i=1}^{r} (d_i)_{y_i} (n_i - y_i)_{x_i},$$

where $y_i$ is the number of neighbors of $u$ that receive color $i$ and $x_i$ is the number of non-neighbors of $u$ that receive color $i$. We find it convenient to work instead with the following polynomial. For $\boldsymbol{\xi} = (\xi_1, \ldots, \xi_r) \in \mathbb{R}^r$, let

$$P_F(\boldsymbol{\xi}) = \frac{1}{\text{Aut}(F)} \sum_{u \text{ critical}} \sum_{\chi_u} \prod_{i=1}^{r} \frac{1}{r^{x_i}} \xi_i^{y_i}.$$

We now state a few easy properties of the polynomial $P_F(\boldsymbol{\xi})$.

**Lemma 2.4.** $P_F(\boldsymbol{\xi})$ *is a symmetric polynomial with non-negative coefficients.* □

**Lemma 2.5.** *For every $\epsilon > 0$, there exists $\delta > 0$ satisfying the following: if $n = \sum_{i=1}^r n_i > 1/\delta$ and if, for all $i \in [r]$, we have $0 \leq d_i \leq n_i$, $|n_i - n/r| \leq \delta n$ and $|\xi_i - d_i/n| \leq \delta$, then $|\#F(n_1, \ldots, n_r; \boldsymbol{d}) - n^{f-1}P_F(\boldsymbol{\xi})| < \epsilon n^{f-1}$.* □

As a first exercise, let us characterize all connected graphs for which $\deg(P_F) = r$ (we will later need to treat such graphs separately).

**Lemma 2.6.** *If $F$ is a connected $r$-critical graph and $\deg(P_F) = r$, then $F = K_{r+1}$ or $r = 2$ and $F = C_{2k+1}$ is an odd cycle.*

*Proof.* The degree of $P_F$ is determined by the largest degree of a critical vertex. Therefore, $\deg(u) \leq r$ for each critical vertex $u \in F$. However, any $r$-coloring $\chi_u$ of $F - u$ must assign all $r$ colors to the neighbors of $u$. Thus, $\deg(u) = r$, every edge incident to $u$ is a critical edge, and, by extension, every neighbor of $u$ is a critical vertex. As $F$ is connected, it follows that every vertex is critical and has degree $r$. Thus the lemma follows from Brooks' Theorem [1]. □

**Lemma 2.7.** *Let $F$ be an $r$-critical graph such that $\deg(P_F) = r$. Then $\pi_F = 1/r$.*

*Proof.* Let $F$ be an $r$-critical graph. Then by Lemma 2.6, we may write $F = F_1 \cup G$ where $F_1$ is a connected $r$-critical graph isomorphic to either $K_{r+1}$ or $C_{2k+1}$ and $G$ is a (possibly empty, not necessarily connected) $r$-colorable graph. Let us first consider the cases where $G$ is empty.

If $F = K_{r+1}$, it is easily seen that $c(n_1, n_2, \ldots, n_r; F) = \prod_{i=2}^r n_i$. By taking $n_i = n/r$, it is immediate that $\alpha_F$, the coefficient of $n^{r-1}$, is $(1/r)^{r-1}$. Next, by taking partial derivatives, we have that $\frac{\partial c}{\partial_2} = \prod_{i=3}^r n_i$ and $\frac{\partial c}{\partial_1} = 0$. Once again letting $n_i = n/r$, we see that $\zeta_F$, the coefficient of $n^{r-2}$ in the difference $\frac{\partial c}{\partial_1} - \frac{\partial c}{\partial_2}$, is $-(1/r)^{r-2}$. Therefore, $\pi_F = \alpha_F/|\zeta_F| = 1/r$.

Next, if $F = C_{2k+1}$, we have $c(n_1, n_2; F) = (n_1 - 2)_{k-1}(n_2)_k$, which is a polynomial of degree $2k - 1$. Thus $\alpha_F = 2^{-(2k-1)}$. Routine calculations show that $\zeta_F = -2^{-(2k-2)}$ and $\pi_F = 1/2 = 1/r$.

Now, if $G$ is not empty, we may write $c(n_1, n_2, \ldots, n_r; F)$ as a product of two polynomials $f$ and $g$, where $f = c(n_1, n_2, \ldots, n_r; F_1)$ and $g$ gives the number of copies of $G$ in the remaining complete $r$-partite graph (that is, copies of $G$ that do not use vertices already claimed by a copy of $F_1$). As calculations of $\alpha_F$, $\zeta_F$ and $\pi_F$ only require the term(s) of highest degree in $c(n_1, n_2, \ldots, n_r; F)$, we will denote by $\hat{c}$, $\hat{f}$, and $\hat{g}$ the respective polynomials consisting of such terms. It follows by our definition that $\hat{c} = \hat{f}\hat{g}$ and, denoting by $d$ the degree of $\hat{g}$ and by $g_0$ the sum of the coefficients of the terms of $\hat{g}$, we have $\alpha_F = g_0(1/r)^d \alpha_{F_1}$.

The polynomial $\hat{g}$ of the highest degree terms of $g$ is symmetric with respect to $n_1, n_2, \ldots, n_r$. Therefore, when we evaluate the polynomials and their derivatives at the vector $(n/r, \ldots, n/r)$, we have that $\frac{\partial \hat{g}}{\partial_1} = \frac{\partial \hat{g}}{\partial_2}$ and

$$\frac{\partial \hat{c}}{\partial_1} - \frac{\partial \hat{c}}{\partial_2} = \hat{f}\left(\frac{\partial \hat{g}}{\partial_1} - \frac{\partial \hat{g}}{\partial_2}\right) + \hat{g}\left(\frac{\partial \hat{f}}{\partial_1} - \frac{\partial \hat{f}}{\partial_2}\right) = \hat{g}\left(\frac{\partial \hat{f}}{\partial_1} - \frac{\partial \hat{f}}{\partial_2}\right).$$

So, $\zeta_F = g_0(1/r)^d \zeta_{F_1}$ and, thus, $\pi_F = \alpha_{F_1}/|\zeta_{F_1}| = \pi_{F_1} = 1/r$. □

**Lemma 2.8.** *Let $F$ be an $r$-critical graph such that $\deg(P_F) = r$. Then $t_F(n, q) = qc(n, F)$ for $q \leq \lfloor n/r \rfloor - 1$.*

*Proof.* Clearly, $t_F(n, q) \geq qc(n, F)$, so we need to construct a graph $H \in \mathcal{T}_r(n, q)$ that has at most $qc(n, F)$ copies of $F$. Take $V(H) = U_1 \cup \ldots \cup U_r$ where $|U_i|$ is either $\lceil n/r \rceil$ or $\lfloor n/r \rfloor$, $c(n, F) = c(|U_1|, \ldots, |U_r|; F)$ and $E(H) = K(U_1, \ldots, U_r) \cup K(\{u\}, W)$, where $u \in U_1$, $W \subseteq U_1 \setminus \{u\}$ and $|W| = q$. That is, $H$ is obtained from $T_r(n)$ by adding (the edges of) a star of size $q$ into $U_1$. Observe that any copy of $F$ in $H$ must use the vertex $u$. Furthermore, $u$ is contained in the $r$-critical component, which is isomorphic to $K_{r+1}$ or $C_{2k+1}$. So, each copy of $F$ uses exactly one bad edge incident to $u$ and $\#F(H) = qc(n, F)$, as required. □

Let us now restrict the domain of $P_F$ to those $\boldsymbol{\xi}$ which may arise as the density vector of some vertex. Note that if $\boldsymbol{d}$ corresponds to the degrees of a vertex and we let $\boldsymbol{\xi} = \boldsymbol{d}/n$, it would follow that $\xi_i \geq 0$ for all $i \in [r]$. Furthermore, as $\sum_{i=1}^{r} d_i \leq n - 1$, we have $\sum_{i=1}^{r} \xi_i \leq 1$. However, we mostly encounter essentially equitable partitions and, therefore, use the more restrictive set

$$\mathcal{S} = \{\boldsymbol{\xi} \in \mathbb{R}^r : \forall i \in [r] \; 0 \leq \xi_i \leq 1/r\}.$$

Many of the arguments that follow involve minimizing $P_F$, usually over some subset of $\mathcal{S}$. One such subset is $\mathcal{S}_\rho = \{\boldsymbol{\xi} \in \mathcal{S} : \sum_{i=1}^{r} \xi_i = \rho + \frac{r-1}{r}\}$ where $\rho \in [0, 1/r]$. Let

$$p(\rho) = \min\{P_F(\boldsymbol{\xi}) : \boldsymbol{\xi} \in \mathcal{S}_\rho\}. \tag{9}$$

**Definition 2.9.** *If $\deg(P_F) \geq r+1$, let $\rho_F = \inf\left\{\rho \in (0, \frac{1}{r}) : p(\rho) \leq \alpha_F \rho\right\}$. If $\deg(P_F) = r$, then define $\rho_F = \infty$.*

**Definition 2.10.** *If $\deg(P_F) \geq r+1$, let $\hat{\rho}_F = \inf\left\{\rho \in (0, \frac{1}{r}) : p(\rho) < \alpha_F \rho\right\}$. If $\deg(P_F) = r$, then define $\hat{\rho}_F = \infty$.*

As it is usual, we agree that the infimum of the empty set is $\infty$. Thus, for example, $\rho_F = \infty$ if $p(\rho) > \alpha_F \rho$ for all $\rho \in (0, \frac{1}{r})$.

Roughly speaking, these definitions consider constructing a graph with $t_r(n) + \rho n$ edges by starting with $T_r(n-1)$ and adding a new vertex $x$ with degrees into the parts being $\xi_1 n + \ldots + \xi_r n = (\rho + \frac{r-1}{r})n$, where $\boldsymbol{\xi} \in \mathcal{S}_\rho$ is optimal, that is, $p(\rho) = P_F(\boldsymbol{\xi})$. The parameter $\rho_F$ (resp., $\hat{\rho}_F$) indicates the first moment when this becomes asymptotically as good as (resp., strictly better than) adding $\rho n$ edges into one part of $T_r(n)$.

**Lemma 2.11.** *If $\deg(P_F) \geq r+1$, then $P_F(\gamma, 1/r, \ldots, 1/r) > \alpha_F \gamma$ for every $\gamma \in (0, 1/r]$.*

*Proof.* One can use the combinatorial interpretation of $P_F$ given by Lemma 2.5. If we connect a vertex $x \in V_1$ in $K(V_1, \ldots, V_r) = T_r(n)$ to $\gamma n + O(1)$ vertices of $V_1 \setminus \{x\}$, then in addition to at least $c(n, F) = \alpha_F n^{f-2} + O(n^{r-3})$ copies of $F$ per each added edge, there will also be $\Omega(n^{f-1})$ copies that use more than one added edge. □

Also, note that $p(0) = 0$. So, for example, Definition 2.10 will not change if the infimum is taken over $\rho \in [0, \frac{1}{r}]$. Clearly, $\rho_F \leq \hat{\rho}_F$. Let us show that $\rho_F$ is strictly positive.

**Lemma 2.12.** *$\rho_F > 0$.*

*Proof.* Assume that $\deg(P_F) \geq r+1$ for otherwise the stated inequality holds by the definition of $\rho_F$. Given $F$, choose sufficiently small positive $\delta$ and then $\epsilon \ll \delta$. Let us show that $\rho_F \geq \epsilon$.

Take any $\rho \in (0, \epsilon)$. It is enough to show that $P_F(\boldsymbol{\xi}) > \alpha_F \rho$ for every $\boldsymbol{\xi} \in \mathcal{S}_\rho$ because then $p(\rho) > \alpha_F \rho$ by the compactness of $\mathcal{S}_\rho$ and the continuity of $P_F$.

As $P_F$ is symmetric, assume without loss of generality that $\xi_1 \leq \xi_i$ for all $i \in [r]$. Also, we may assume that $\xi_1 < \delta$ for otherwise, since $P_F$ has non-negative coefficients, we are done:

$$P_F(\boldsymbol{\xi}) \geq P_F(\delta, \ldots, \delta) > \epsilon \alpha_F > \rho \alpha_F.$$

Fix some index $i$ with $2 \leq i \leq r$. Since $\sum_{j=2}^{r}(\frac{1}{r} - \xi_j) = \xi_1 - \rho < \delta$, we have that $0 \leq 1/r - \xi_i < \delta$. Since each monomial of $P_F(\boldsymbol{\xi})$ contains $\xi_1$, the $i$-th partial derivative of $P_F$ at $\boldsymbol{x_0} = (0, \frac{1}{r}, \ldots, \frac{1}{r})$ is 0. Thus, as $\frac{\partial}{\partial_1} P_F(\boldsymbol{x_0}) = \alpha_F > 0$ and $\delta$ is small, we have that

$$\frac{\partial}{\partial_1} P_F(\boldsymbol{\xi}) > \frac{\partial}{\partial_i} P_F(\boldsymbol{\xi}) + \frac{\alpha_F}{2}.$$

If follows that, if we increase $\xi_i$ to $1/r$ and decrease $\xi_1$ by the same amount, then the value of $P_F(\boldsymbol{\xi})$ does not go up. Iteratively repeating the above perturbation for each $i \geq 2$, we obtain that $P_F(\boldsymbol{\xi}) \geq P_F\left(\rho, \frac{1}{r}, \ldots, \frac{1}{r}\right)$ which is strictly larger than $\rho \alpha_F$ by Lemma 2.11, as required. □

To give a better picture of proceedings, let us informally recall some previous parameters. (Also, the glossary at the end of the paper can be used for looking up many definitions.) One of our main objectives is to estimate the parameters $c_1(F)$ and $c_2(F)$ that were defined respectively in (4) and (3). First, consider starting with the Turán graph and growing the graph by adding extra edges optimally. When scaled appropriately, the number of copies of $F$ grows 'linearly' in $q$ with a slope of $\alpha_F > 0$ (defined in Lemma 2.1). On the other hand, if we start with a slight perturbation of the partition sizes, we may have a slope slightly smaller than $\alpha_F$ but a higher intercept. The ratio $\pi_F$ (Definition 2.3) gives the intersection of these two 'ideal' lines, where $F$-subgraphs that use more than one extra edge are ignored. Alternatively, we may start with a Turán graph on one fewer vertices and grow the graph by introducing a vertex of appropriate degree. The number of copies then grows according to $p(\frac{q}{n})$. In this scenario, $\rho_F$ and $\hat{\rho}_F$ (Definitions 2.9–2.10) identify the first time when this curve, respectively, intersects and crosses the line of slope $\alpha_F$. (Thus $c_2(F) \leq \hat{\rho}$ and our Theorem 3.12 shows that this is equality.)

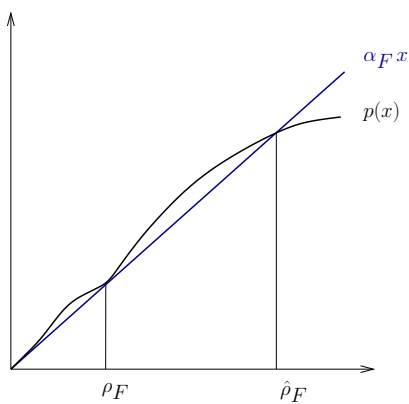

Figure 1: $\rho_F$ and $\hat{\rho}_F$.

The values $\rho_F$ and $\hat{\rho}_F$ in Figure 1 do not coincide. However, $\rho_F = \hat{\rho}_F$ for all graphs we have thus far encountered, and it may be possible that equality holds for all graphs. In some instances, this would imply that $c_1(F) = c_2(F)$.

## 3. General Results

Here we prove Theorem 1.4 with an explicit lower bound on $c_1(F)$. Then we show that, for $q = o(n^2)$, the limit of $h_F(n,q)/qn^{f-2}$ can be computed as the value of some optimization problem, see Theorems 3.10 and 3.11. This will allow us to determine $c_2(F)$ in Theorem 3.12. Our approach is based on Lemmas 3.3 and 3.4 on the structure of $F$-optimal graphs. As was the case in Mubayi's paper [16], the Graph Removal Lemma (see e.g. [11, Theorem 2.9]) and the Erdős-Simonovits Stability Theorem are key components of our proof.

**Theorem 3.1** (Graph Removal Lemma). *Let $F$ be a graph with $f$ vertices. Then for every $\epsilon_1 > 0$ there is $\epsilon_2 > 0$ such that every graph $H$ with $n \geq 1/\epsilon_2$ vertices and at most $\epsilon_2 n^f$ copies of $F$ can be made $F$-free by removing at most $\epsilon_1 n^2$ edges.*

**Theorem 3.2** (Erdős [6] and Simonovits [21]). *Let $r \geq 2$ and $F$ be a graph with chromatic number $r + 1$. Then for every $\epsilon_1 > 0$ there is $\epsilon_2 > 0$ such that every $F$-free graph $H$ with $n \geq 1/\epsilon_2$ vertices and at least $t_r(n) - \epsilon_2 n^2$ edges contains an $r$-partite subgraph with at least $t_r(n) - \epsilon_1 n^2$ edges.*

### 3.1. Common Definitions and Auxiliary Results

Here we provide some definitions and auxiliary results that will be used by Theorems 3.5–3.10 and also in Section 4.

Let us start by defining some constants. Given an $r$-critical graph $F$, pick constants satisfying the following hierarchy:

$$\delta_0 \gg \delta_1 \gg \delta_2 \gg \delta_3 \gg \delta_4 \gg \delta_5 \gg \delta_6 \gg \delta_7 \gg \delta_8 \gg 1/n_0, \tag{10}$$

each being sufficiently small depending on the previous ones. While most of the dependencies can be explicitly written, one of them depends on the Removal Lemma. Thus, the final $n_0 = n_0(F)$ is very large and there is not much point in optimizing these constants.

The following definitions will apply to the rest of the paper, given an optimal graph $H \in \mathcal{H}_F(n,q)$ where $n \geq n_0$ and $1 \leq q < \delta_8 n^2$. We fix a *max-cut $r$-partition* $V = V_1 \cup \ldots \cup V_r$ of the vertex set $V = V(H)$, that is, we maximize $|H \cap T|$, where $T = K(V_1, \ldots, V_r)$ is the complete $r$-partite graph with parts $V_i$. We let $B = H \setminus T$ and call the edges in $B$ *bad*. Thus, an edge of $H$ is bad if it lies inside some part $V_i$. Let $M = T \setminus H$ be the set of *missing* edges. Since $|T| \leq t_r(n)$, we have that

$$|B| \geq |M| + q. \tag{11}$$

Note that every copy of $F$ in $H$ must contain at least one bad edge. To this end, we denote by $\#F(uv)$ the number of $F$-subgraphs of $H$ that contain the bad edge $uv \in B$ but no other bad edges. In addition, for $u \in V(H)$, the number of copies of $F$ in $H$ that use the vertex $u$ is denoted by $\#F(u)$.

If we consider a graph in $\mathcal{T}_r(n,q)$ so that the maximum bad degree is at most, say, $2q/n + 2$, then it has at most $2c(n;F)q + 2f^4q^2n^{f-4} \leq 3\alpha_F \delta_8 n^f$ copies of $F$. (For example, the number of $F$-subgraphs that contain two disjoint bad edges is at most $\binom{q}{2}f^4 n^{f-4}$.) It follows that $h_F(n,q) \leq t_F(n,q) \leq 3\alpha_F \delta_8 n^f$. Thus the Removal Lemma applies to $H$ and gives an $F$-free subgraph $H' \subseteq H$ with at least $t_r(n) - \delta_7 n^2$ edges. We then apply the Erdős-Simonovits Stability Theorem and obtain an $r$-partite subgraph $H'' \subseteq H'$ with at least $t_r(n) - \delta_6 n^2$ edges. By the max-cut property, we have that $|H \cap T| \geq |H''|$. It routinely follows that

$$|V_i| = \left(\frac{1}{r} \pm \delta_5\right) n. \tag{12}$$

Also, we have

$$|B| \leq q + \delta_6 n^2 \leq 2\delta_6 n^2. \tag{13}$$

For $uv \in M$, define

$$\#F'(uv) = \#F(H + uv) - \#F(H)$$

to be the number of *potential* copies of $F$ associated with $uv$, that is, the number of copies of $F$ introduced by including the missing edge $uv$. Note that

$$\#F(xy) \leq \#F'(uv), \qquad \text{for all } uv \in M, \ xy \in B, \tag{14}$$

as otherwise we can reduce the number of copies of $F$ in $H$ by removing the edge $xy$ and adding $uv$.

Let $d$ denote the degree of a vertex; for example, $d_B(x)$ is the number of bad edges containing $x$.

**Lemma 3.3.** *If $M \neq \emptyset$, then there exists a vertex $x \in V$ with $d_B(x) \geq \delta_1 n$.*

*Proof.* Assume, for contradiction, that $d_B(u) < \delta_1 n$ for all $u \in V$. If $uv \in M$, then it follows from $d_B(u) + d_B(v) < 2\delta_1 n$ that $\#F'(uv) < 2\delta_1 f^3 n^{f-2} + qf^4 n^{f-4}$. Let $xy \in B$ be arbitrary. Since $|M| \leq |B| \leq 2\delta_6 n^2$ by (11) and (13), we have that

$$\#F(xy) \geq \alpha_F n^{f-2} - (d_M(x) + d_M(y)) f^3 n^{f-3} - 2\delta_4 f^4 n^{f-2}. \tag{15}$$

Now, (14) implies that $d_M(x) + d_M(y) \geq (\alpha_F - \delta_0) n / f^3$. By counting the number of adjacent pairs $(e, e')$ with $e \in B$ and $e' \in M$ in two different ways, we obtain

$$|M| \times 2\Delta(B) \geq |B| \times \frac{(\alpha_F - \delta_0)n}{f^3}.$$

Since $\Delta(B) < \delta_1 n$, we conclude that $|B| \leq |M|$, contradicting (11). □

**Lemma 3.4.** *Let $x$ be a vertex with $d_B(x) \geq \delta_1 n$. Then $\#F(x) \geq \delta_2 n^{f-1}$ and $d_H(x) \geq \left(\frac{r-1}{r} + \delta_3\right) n$.*

*Proof.* As $V_1 \cup \ldots \cup V_r$ is a max-cut partition, we have $|N(x) \cap V_i| \geq \delta_1 n$ for all $i \in [r]$. It follows from Lemma 2.5 that
$$\#F(x) \geq P_F(\delta_1, \ldots, \delta_1) n^{f-1} - \delta_2 n^{f-1} \geq \delta_2 n^{f-1}.$$

Furthermore, as $P_F(\boldsymbol{\xi})$ is continuous and $P_F(0, 1/r, \ldots, 1/r) = 0$, we have that $P_F(\delta_3 + (r-1)\delta_5, 1/r, \ldots, 1/r) \leq \delta_2/2$. Hence, unless $d_H(x) \geq (\frac{r-1}{r} + \delta_3)n$, we can strictly reduce the number of copies of $F$ by connecting $x$ to all vertices in other classes and to at most $(\delta_3 + (r-1)\delta_5)n$ vertices in its own class. □

3.2. *Lower Bound on $c_1(F)$*

Lemmas 3.3 and 3.4 directly imply that, if $q < \delta_3 n$, then no $H \in \mathcal{H}_F(n, q)$ can have a missing edge. This is a big step towards proving Theorem 1.4. However, we aim for a better explicit bound on $c_1(F)$, which requires considering also the case $M \neq \emptyset$ in the proof. Recall that $c_{1,i}(F)$ was defined in (5). Let $\mathrm{sign}(x)$ denote the sign of $x \in \mathbb{R}$: it is 0 if $x$ is 0 and $x/|x|$ otherwise.

**Theorem 3.5.** *Let $F$ be an $r$-critical graph and let $t = -\mathrm{sign}(\zeta_F)$. Then the following claims hold.*

- *If $\deg(P_F) = r$, then $c_1(F) \geq 1/r$.*
- *If $\deg(P_F) \geq r+1$, then $c_1(F) \geq \min(\pi_F, \rho_F)$ and $c_{1,t}(F) \geq \min(2\pi_F, \rho_F)$.*

*Proof.* As stated in the Introduction, we prove the theorem by showing something stronger; namely, for all large $n$ with $q$ in the appropriate range, not only is $h_F(n, q) = t_F(n, q)$ but also, in fact, $\mathcal{H}(n, q) \subseteq \mathcal{T}_r(n, q)$. In particular, this will imply Theorem 1.4 as $\rho_F > 0$ by Lemma 2.12.

Let $\epsilon > 0$ be arbitrary. Consider $H \in \mathcal{H}_F(n, q)$ where $n$ is large and $q < n/\epsilon$. Let all definitions of Section 3.1 apply to $H$, where we require that $\delta_0 \ll \epsilon$. Also, let $c(n, F) = c(n'_1, \ldots, n'_r; F)$ be as in (7).

Consider the following cases:

**Case 1:** $M = \emptyset$.

Thus we know that $H \supseteq K(V_1, \ldots, V_r)$ and our aim is to conclude that the part sizes are nearly equal. By symmetry, assume that $|V_1| \geq \cdots \geq |V_r|$. We show that if $|V_1| \geq |V_r| + 2$ then $q \geq (c - \epsilon)n$ where $c$ is the lower bound on $c_1(F)$ that we have to prove. In other words, if $q < (c - \epsilon)n$, then $H \in \mathcal{T}_r(n, q)$.

Let $a = \max(|V_1| - \lceil n/r \rceil, \lfloor n/r \rfloor - |V_r|)$. We have
$$q + \frac{a^2 r}{2(r-1)} \leq |B| \leq q + \frac{(a^2 + a)r}{2}.$$

As $a \leq \delta_5 n$ by (12), Lemma 2.2 implies that, regardless of how the bad edges are distributed among parts, we have
$$\begin{aligned}
\#F(H) &\geq \left(q + \frac{a^2 r}{2(r-1)}\right)\left(c(n, F) - |\zeta_F| a n^{f-3} - \frac{a^2 n^{f-4}}{\delta_0}\right) \\
&\geq qc(n, F) + an^{f-2}\left(\alpha_F \frac{ar}{2(r-1)} - \frac{q}{n}|\zeta_F| - \delta_1 a\right).
\end{aligned}$$

On the other hand, $t_F(n, q) \leq qc(n, F) + n^{f-2}/\delta_0$, which is demonstrated by adding extra $q$ edges to one part of $T_r(n)$ so that they form a graph of bounded maximum degree. Since $H$ is optimal, it follows that $a < 1/\delta_1$. In particular, we have $|B| \leq q + 1/\delta_2$.

For $i \in [r]$, let $B_i = B[V_i]$ be the set of bad edges contained in $V_i$.

**Claim 3.6.** *If $|V_j| = |V_k| + s$, where $s > 1$, then*

$$(|B_k| - |B_j|)\zeta_F \geq (1 - \delta_2)(s-1)\alpha_F n. \tag{16}$$

*Proof of Claim.* Consider moving one vertex $v$ from $V_j$ to $V_k$ as follows.

Pick a vertex $v \in V_j$ with $d_B(v) \leq d_B(u)$ for all $u \in V_j$. Thus $d_B(v) = O(1)$. Move all $\ell = d_B(v)$ bad edges incident to $v$ inside $\{x \in V_j \setminus \{v\} : d_B(x) \leq 1/\delta_0\}$. Since $q/n < 1/\epsilon$, this set contains at least, say, $n/2r$ vertices and this move is possible. The number of $F$-subgraphs increases by at most $\ell f^4(\delta_0^{-1} n^{f-3} + qn^{f-4}) = O(n^{f-3})$. Let us update $B_i$'s to refer to the new graph. Note that the sizes $|B_i|$ do not change.

Next, we replace $v$, which is an isolated vertex in the graph $B$ now, with a vertex $v'$ such that $uv' \in H'$ for all $u \in V \setminus V_k$. By Lemma 2.2, the $F$-count changes by $(|B_k| - |B_j|)\zeta_F n^{f-3} + O(n^{f-3})$. (For example, $\#F(xy)$ changes by at most $O(n^{f-4})$ for every bad edge $xy$ not inside $V_j \cup V_k$.)

Finally, as $(|V_j| - 1)(|V_k| + 1) = |V_j||V_k| + s - 1$, we remove $s - 1$ bad edges chosen arbitrarily. By (12), the $F$-count drops by at least $(1 - \delta_3)(s-1)c(n, F)$.

Since we have arrived at a graph with the same edge count as $H \in \mathcal{H}_F(n, r)$, the number of $F$-subgraphs cannot decrease during the above transformations. Thus

$$(|B_k| - |B_j|)\zeta_F n^{f-3} - (1 - \delta_3)(s-1)c(n, F) \geq O(n^{f-3}),$$

which implies the claim. ∎

The above claim and our assumption $|V_1| \geq |V_r| + 2$ imply that $|B| \geq \max(|B_1|, |B_n|) \geq (1 - \delta_2)\pi_F n$, and $q \geq |B| - (a^2 + a)r/2 \geq (\pi_F - \delta_0)n$, as desired. (Recall that if $\deg(P_F) = r$, then $\pi_F = 1/r$ by Lemma 2.7.) Also, it follows that $\zeta_F \neq 0$. Thus it remains to show that $c_{1,t}(F) \geq 2\pi_F$, where $t = -\operatorname{sign}(\zeta_F) \neq 0$. We modify the previous argument as follows. Let $n \equiv t \pmod{r}$. Then there exist $j, k, l$ with $j \neq k$ such that $|V_j| - |V_k| \geq 3$ or $|V_j| - |V_l| = |V_k| - |V_l| = 2t$. In the first case, we apply Claim 3.6 as above to obtain $q \geq (2\pi_F - \delta_0)n$. In the second case, we have, by applying Claim 3.6 twice, that $|B_j|, |B_k| \geq (1 - \delta_2)\pi_F n$, again implying that $q \geq (2\pi_F - \delta_0)n$.

**Case 2:** $M \neq \emptyset$.

It follows from Lemma 3.3 that $X \neq \emptyset$, where we define

$$X = \{x \in V(H) : d_B(x) \geq \delta_1 n\}. \tag{17}$$

We will now handle the two cases $\deg(P_F) = r$ and $\deg(P_F) \geq r + 1$ separately. Let us first consider the case $\deg(P_F) \geq r + 1$.

**Claim 3.7.** *If $\deg(P_F) \geq r + 1$, then $d_H(x) \geq (\rho_F + \frac{r-1}{r} - \delta_2)n$ for all $x \in X$. In particular, since $X \neq \emptyset$, we have that $\rho_F \neq \infty$ (and thus $\rho_F < 1/r$).*

*Proof of Claim.* Let, for example, $x \in X \cap V_1$ contradict the claim. By the definition of $\rho_F$, we have $p(\rho) > \alpha_F \rho$ for $0 < \rho < \min(\rho_F, 1/r)$, where $p$ is defined by (9). Since $p$ is a continuous function, we can assume that $p(\gamma) - \alpha_F \gamma \geq 5\delta_3$, where we let $\gamma = \frac{d_H(x)}{n} - \frac{r-1}{r}$. (Note that (12) and Lemma 3.4 show that $\gamma \geq \delta_3 - (r-1)\delta_5$ is separated from 0 while e.g. Lemma 2.11 takes care of the case when $\gamma$ is close to $1/r$.)

Let us replace $x$ with a vertex $u$ whose neighborhood is $V(H) \setminus V_1$. Clearly, the new count of $F$-subgraphs satisfies $\#F(u) \leq |B|f^3 n^{f-3} \leq \delta_3 n^{f-1}$. Next, we distribute the remaining $d(x) - d(u)$ edges evenly among vertices in $V_1$ with bad degree at most $\delta_5 n$. This creates at most $(\gamma n + r\delta_5 n)c(n, F) + 2\delta_3 n^{f-1}$ copies of $F$. Comparing with the old value $\#F(x) \geq (p(\gamma) - \delta_3)n^{f-1}$, we see that the number of $F$-subgraphs decreases, a contradiction to the optimality of $H$. ∎

Note by Claim 3.7 that if $x \in X$ then

$$\#F(x) \geq n^{f-1}(p(\rho_F - 2\delta_2) - \delta_3) > n^{f-1}(\alpha_F \rho_F - \delta_1) > (\rho_F - \delta_0)nc(n, F) + \delta_1 n^{f-1}. \tag{18}$$

Thus, if $q \leq (\rho_F - \delta_0)n$, the number of copies of $F$ at some vertex $x \in X$ alone exceeds the upper bound of $qc(n, F) + O(qn^{f-3})$ on $h_F(n, q)$, contradicting our assumption that $\#F(H) = h_F(n, q)$. This completes the proof of Theorem 3.5 for $r$-critical graphs with $\deg(P_F) \geq r + 1$.

We now consider the case when $\deg(P_F) = r$. Here we take $H \in \mathcal{H}_F(n, q)$ where $q \leq (1/r - \delta_0)n$. Recall the set $X$ defined in (17).

**Claim 3.8.** *Every $x \in X$ is incident to at most $(r - 1)\delta_2 n$ missing edges.*

*Proof of Claim.* Let $x \in X \cap V_1$ with $d_i$ neighbors in each part $V_i$. By the max-cut property we have $d_1 \leq d_i$ for all $i \in [r]$. In particular, each $d_i$ is at least $\delta_1 n$.

If, for example, $x$ has at least $\delta_2 n$ non-neighbors in $V_2$, then we can move $\delta_2 n/2$ edges at $x$ from $V_1$ to $V_2$, decreasing the product $d_1 d_2$ by at least $\delta_3 n^2$. Since $P_F(\boldsymbol{\xi}) = C_F \prod_{i=1}^{r} \xi_i$ for some constant $C_F > 0$, this would strictly decrease $\#F(x)$, a contradiction to the extremality of $H$. The claim follows. ∎

**Claim 3.9.** *Every missing edge intersects $X$.*

*Proof of Claim.* Suppose on the contrary that there exists $uv \in M$ with $u, v \notin X$. As both endpoints have bad degree of at most $\delta_1 n$, it follows that $\#F'(uv) \leq 2\delta_1 f^3 n^{f-2} + qf^4 n^{f-4}$. On the other hand, consider a vertex $x \in X$. There is a bad edge $xw$ such that $d_M(w) < \delta_3 n$ for otherwise we would have $|B| \geq |M| \geq \frac{1}{2} d_B(x)\delta_3 n > 2\delta_6 n^2$, contradicting (13). But then, by Claim 3.8, we have $\#F(xw) \geq \alpha_F n^{f-2} - \delta_1 n^{f-2} > \#F'(uv)$, contradicting (14). ∎

As $\#F(H) \leq 3\alpha_F \delta_8 n^f$, it follows from Lemma 3.4 that $|X| < \delta_5 n$. Thus, by Claims 3.8 and 3.9 we have that
$$\Delta(M) \leq \max\Big(|X|, (r-1)\delta_2 n\Big) = (r-1)\delta_2 n.$$

It follows that $\#F(u'v') \geq (1 - \delta_1)c(n, F)$ for every bad edge $u'v' \in B$. Since $\#F(H) \leq (1 + \delta_1)qc(n, F)$, we conclude that
$$|B| \leq \frac{(1 + \delta_1)q}{1 - \delta_1} < (1/r - \delta_0/2)n. \tag{19}$$

Now pick some vertex $u^*$ and consider a graph $H'$ where all $|B|$ bad edges and $|M|$ missing pairs are incident to $u^*$. This eliminates all copies of $F$ using multiple bad edges. Furthermore, observe that if a missing pair $uv$ and a bad edge $u'v'$ in $H$ are disjoint, then we increased the number of potential copies that contain both $uv$ and $v'v'$ by $\Omega(n^{f-3})$ when we made them adjacent.

Therefore, by choosing the part containing $u^*$ appropriately, we achieve that $\#F_H(u'v') \geq \#F_{H'}(u^*w)$ for every $u'v' \in B(H)$ and $u^*w \in B(H')$. Thus $\#F(H') \leq \#F(H)$. This has to be equality by the extremality of $H$. So we can assume that $H = H'$, that is, $H$ is obtained by adding a vertex $u^*$ to a complete $r$-partite graph of order $n - 1$. Since $p(x)$ is a continuous function, we can assume by Lemma 2.7 that $p(x) \geq \alpha_F x + \delta_4$ for every $x \in [\delta_3, \frac{1}{r} - \delta_3]$. Now, it follows from Lemmas 2.5 and 3.4 and from (12) that the degree of $u^*$ should be at least $(1 - \delta_1)n$.

This implies that $|B| \geq (1/r - \delta_1)n$, which contradicts our assumption on $q$ and completes the proof of Theorem 3.5. □

### 3.3. Vertex Removal Procedure

In this section, we describe a procedure that, given a graph $H$, produces a sequence of numbers that can be used to lower bound the number of $F$-subgraphs in $H$. Conversely, one can 'assemble' a graph in an almost optimal way for a given sequence. This allows us to write the limit of $h_F(n, q)/qn^{f-1}$ for $q = o(n^2)$ as the value of a certain optimization problem.

Fix constants as in (10); they will not be modified during the procedure. Let $n \geq 2n_0$, $q \leq \delta_8 n^2/4$, and $H \in \mathcal{H}_F(n, q)$ be arbitrary.

**Vertex Removal Procedure.** Initially, let $H_0 = H$ and $q_0 = q$. Iteratively for $i = 1, 2, \ldots$, given a graph $H_i \in \mathcal{H}_F(n - i, q_i)$, we either stop or construct $H_{i+1}$ using the following steps:

1. If $i \geq n/2$ or $q_i < 0$, then we let $k = i$ and stop.
2. Let all definitions and results of Section 3.1 apply to $H_i \in \mathcal{H}_F(n-i, q_i)$. (Note that $|V(H_i)| = n - i > n/2 \geq n_0$ while Inequality (20) implies that $q_i \leq q \leq \delta_8 n^2/4 \leq \delta_8(n-i)^2$.)
3. If there are no missing edges, let $k = i$ and stop.
4. Pick an arbitrary vertex $x_i \in V(H_i)$ with bad degree at least $\delta_1(n-i)$. (It exists by Lemma 3.3.)
5. Let $q_{i+1} = |H_i - x_i| - t_r(n-i)$ and take an arbitrary $H_{i+1} \in \mathcal{H}_F(n-i-1, q_{i+1})$, i.e., $H_{i+1}$ is a graph of the same order and size as $H_i - x_i$ that minimizes the number of $F$-subgraphs. (If $q_{i+1} \leq 0$, then $H_{i+1}$ is $F$-free.)
6. Increase $i$ by 1 and repeat the above steps.

The obtained numbers $q_0 = q, q_1, \ldots, q_k$, even though they may depend on the choices made during the procedure, will be very useful to us.

Combining the trivial upper bound $d_{H_i}(x_i) \leq n - i - 1$ and the lower bound of Lemma 3.4, we conclude that, for every $0 \leq i < k$,

$$\delta_3(n-i) \leq q_i - q_{i+1} \leq (n-i-1) - t_r(n-i) + t_r(n-i-1) = \frac{n-i}{r} \pm 1. \tag{20}$$

Let us denote $\tau_i = (q_i - q_{i+1})/(n-i)$ for $0 \leq i < k$ and $\tau_k = \max(0, q_k/(n-k))$. We have that

$$q \leq \sum_{i=0}^{k} \tau_i(n-i) \tag{21}$$

Since the numbers $q_i$ decrease at rate at least $\delta_3(n-k)$ by (20) and $q_{k-1} \geq 0$, we conclude that

$$k \leq \frac{2q}{\delta_3 n} \leq \delta_7 n. \tag{22}$$

In particular, $k < n/2$. Thus the procedure stops because the final $H_k \supseteq K(V_1, \ldots, V_r)$ has no missing edges or $q_k < 0$. It follows that

$$\#F(H_k) \geq (1 - \delta_5)\alpha_F \tau_k n^{f-1}.$$

For $0 \leq i < k$, there are at least $(p(\tau_i) - \delta_4)n^{f-1}$ copies of $F$ in $H_i$ that contain $x_i$. (Note that $d_{H_i}(x_i) = (\tau_i + \frac{r-1}{r})(n-i) + O(1)$; also, by (11) and (13), there are at most $2\delta_6 n^2$ missing edges in $H_i - x_i$ and each can destroy at most $f^3 n^{f-3}$ copies of $F$ via $x_i$.) Hence, the number of $F$-subgraphs in the initial graph $H_0 = H$ is

$$\#F(H) \geq \alpha_F \tau_k n^{f-2} + \sum_{i=0}^{k-1} p(\tau_i)n^{f-1} - \delta_3 q n^{f-2}. \tag{23}$$

On the other hand, given numbers $q_0, q_1, \ldots, q_k$ that satisfy $q_0 \leq \delta_8 n^2/4$, $q_{k-1} \geq 0$ and (20), one can construct a graph with $n$ vertices and at least $t_r(n) + q_0$ edges as follows.

We start with the Turán graph $K(V_1, \ldots, V_r) \cong T_r(n-k)$. If $q_k > 0$, we add $q_k$ extra edges so that they form a graph of maximum degree at most $2q_k/n + 2$, creating at most $(\alpha_F + \delta_7)q_k n^{f-2}$ copies of $F$.

Then, iteratively for $i = k-1, \ldots, 0$, let

$$d_i = (t_r(n-i) + q_i) - (t_r(n-i-1) + q_{i+1}),$$

add a new vertex $x_i$ and add $\min(d_i, n-k)$ edges between $x_i$ and $V_1 \cup \cdots \cup V_r$ so that $K(V_1, \ldots, V_r) + x_i$ has the smallest number of $F$-subgraphs via $x$. This number is at most $(p(\tau_i) + \delta_2)n^{f-1}$, where $\tau_i$ is defined as before (21).

Let $H'$ be the obtained graph, after we added $x_{k-1}, \ldots, x_0$. Finally, let $H$ be obtained by adding $q' = \max(0, t_r(n) + q_0 - |H'|)$ arbitrary edges to $H'$ (so that $H$ has at least $t_r(n) + q$ edges).

Let us show that $q'$ is small relative to $q$. Note that $q'$ is at most the total number of 'surplus' edges over those vertices $x_i$ for which $d_i > n - k$. Recall that $d_i \le n - i - 1$ by (20). Also, if $d_i > n - k$, then we have e.g. $q_{i+1} \le q_i - n/2r$ and thus, by (22),

$$\frac{d_i - n + k}{q_i - q_{i+1}} \le \frac{k-1}{q_i - q_{i+1}} \le 2r\delta_7.$$

Since the sequence $(q_i)_{i=0}^k$ is monotone decreasing, we have that $q' \le 2r\delta_7 q_0$.

Therefore, the constructed graph $H$ of size $t_r(n) + q$ satisfies

$$\#F(H) \le \alpha_F \tau_k n^{f-1} + \sum_{i=0}^{k-1} p(\tau_i) n^{f-1} + 2r\delta_7 q_0 f^2 n^{f-2} + k\delta_2 n^{f-1}. \tag{24}$$

The above discussion allows us to determine $h_F(n,q)$ asymptotically for $q = o(n^2)$, modulo solving some numerical optimization problem. With this in mind, define $\beta_F$ to be the infimum of the ratio $p(x)/x$ over $x \in (0, \frac{1}{r})$. Observe that $p(x) = \alpha_F x + O(x^2)$ for $x \to 0$. Thus $\beta_F = \alpha_F$ if $\hat{\rho}_F = \infty$ (and $\beta_F < \alpha_F$ otherwise).

**Theorem 3.10.** *If $q = o(n^2)$ and $q/n \to \infty$, then*

$$h_F(n,q) = (\beta_F + o(1))q n^{f-2}.$$

*Proof of Theorem 3.10* Given $q = q(n) = o(n^2)$ and $\epsilon > 0$, choose sufficiently small constants as in (10). Let $n \ge 2n_0$.

To show the lower bound, take an arbitrary $H \in \mathcal{H}_F(n,q)$. Apply the Vertex Removal Procedure to $H$. The definition of $\beta_F$ implies that that $p(x) \ge \beta_F x$ for all $x \in (0, \frac{1}{r})$. Now, (21), (22) and (23) imply that $\#F(H) \ge (\beta_F - \epsilon) q n^{f-2}$, giving the required.

Roughly speaking, to show the upper bound, we fix $\tau \in (0, 1/r)$ such that $p(\tau)/\tau$ is close to its infimum $\beta_H$ and, starting with $q_0 = q$, iteratively decrease $q_i$ in steps $\tau n$ as long as possible. Specifically, choose $\tau \in (2\delta_3, \frac{1}{r} - \delta_3)$ with $p(\tau) \le \beta_F \tau + \epsilon/3$, which is possible since $\delta_3 \ll \epsilon$. Starting with $q_0 = q$, define inductively $q_{i+1} = q_i - \lfloor \tau(n-i) \rfloor$ unless this makes $q_{i+1}$ negative when we stop and let $k = i$. This sequence satisfies (20). Construct $H$ as above. By (24) and since $q_k = O(n)$, we conclude that $h_F(n,q) \le (\beta + \epsilon) q n^{f-2}$, finishing the proof. $\square$

By Theorems 3.5 and 3.10, it remains to consider the case $n/C \le q \le Cn$ for some constant $C > 0$. Then $k \le C/\delta_3$ and determining the asymptotics of $h_F(n,q)$ reduces to the following optimization problem. Let $c = q/n$ and define

$$\phi(c) = \min_k \min_{\boldsymbol{\tau}} \left( \tau_k \alpha_F + \sum_{i=0}^{k-1} p(\tau_i) \right), \tag{25}$$

with the minimum taken over all integers $0 \le k \le C/\delta_3$ and all vectors $\boldsymbol{\tau} \in \mathbb{R}^{k+1}$ satisfying $\sum_{i=0}^k \tau_i \ge c$, $\tau_k \ge 0$, and $\tau_i \ge \delta_3$ for $0 \le i \le k-1$. Then the following holds.

**Theorem 3.11.** *For any constants $C, \epsilon > 0$ and a function $q = q(n)$ such that $n/C \le q \le Cn$ we have that $h_F(n,q) = (\phi(q/n) \pm \epsilon) n^{f-1}$ for all large $n$.* $\square$

3.4. *Determination of $c_2(F)$*

It is not hard to show that the function $\phi$ defined in (25) is continuous. (This also follows from Theorem 3.11.) Thus $c_2(F)$ is the supremum of $c > 0$ for which $\phi(c) = \alpha_F c$. We can in fact pinpoint this value exactly:

**Theorem 3.12.** *For every $r$-critical $F$, we have that $c_2(F) = \hat{\rho}_F$. Furthermore, if $\hat{\rho}_F = \infty$, then $h_F(n,q) = (1 + o(1)) q c(n, F)$ for all $q = o(n^2)$.*

*Proof.* Let us first show that $c_2(F) \leq \hat{\rho}_F$. We may assume that $\hat{\rho}_F$ is finite, for otherwise the upper bound holds vacuously. Let $c > \hat{\rho}_F$ be arbitrary. Take $\boldsymbol{\xi} \in \mathcal{S}_\rho$ such that $\hat{\rho}_F < \rho < \min(c, 1/r)$ and $\lambda > 0$, where $\lambda = \alpha_F \rho - P_F(\boldsymbol{\xi})$.

Let $n$ be large. Let $H$ be obtained from $T_r(n-1)$ by adding a new vertex $u$ that has $(\xi_i + o(1))n$ neighbors in each part $V_i$. Thus $H$ has $t_r(n) + q$ edges, where $q = (\rho + o(1))n$. Then

$$\#F(H) = (P_F(\boldsymbol{\xi}) + o(1))\, n^{f-1} < (\alpha_F \rho - \lambda/2)\, n^{f-1} < (1 - \lambda/3) qc(n, F),$$

This infinite sequence of graphs implies the stated upper bound on $c_2(F)$.

Conversely, let $\epsilon > 0$ be arbitrary and take any function $q = q(n)$ such that $q = o(n^2)$ if $\hat{\rho}_F = \infty$ and $q \leq (\hat{\rho}_F - \epsilon)n$ otherwise. Let $n$ be large and $H \in \mathcal{H}_F(n, q)$. We have to show that $\#F(H) \geq (1 - \epsilon) qc(n, F)$.

Apply the Vertex Removal Procedure to $H$ (where we assume that $\delta_0 \ll \epsilon$). First, if $k = 0$ (in other words, if $H$ has no missing edges), then $\#F(H) \geq (1 - \epsilon) qc(n, F)$ by (12), and we are done. So, suppose now that $k \geq 1$. If there exists some $i$ with $d_{H_i}(x_i) \geq (\hat{\rho}_F + (1 - 1/r))(n - i)$, then, by monotonicity of $p(\rho)$, we have that

$$\#F(H) \geq \#F(x_i; H_i) \geq (p(\hat{\rho}_F) - \delta_3)(n-i)^{f-1} > (1-\epsilon) qc(n, F),$$

that is, the vertex $x_i$ alone provides the required number of $F$-subgraphs. Finally, if $d_{H_i}(x_i) < (\hat{\rho}_F + (1 - 1/r))(n - i)$ for all $i \leq k$, then

$$\#F(x_i; H_i) \geq \alpha_F(d_{H_i}(x_i) - (1 - 1/r)n) n^{f-2} - \delta_3 n^{f-1}.$$

We get the required inequality by summing this quantity over all vertices $x_i$ as in (23). $\square$

## 4. Special Graphs

In this section we obtain upper bounds on $c_{1,i}(F)$ for a class of graphs and compute the exact value for some special instances. We also give an example of a graph with $c_1(F)$ strictly greater than $\min(\pi_F, \rho_F)$.

### 4.1. $K_{r+2} - e$.

Let $r \geq 2$ and let $F = K_{r+2} - e$ be obtained from the complete graph $K_{r+2}$ by deleting one edge. Here we prove Theorem 1.6. Namely, we show that $c_{1,t}(F) = \frac{r-1}{r^2}$ if $t \not\equiv 1 \pmod{r}$ and $c_{1,1}(F) = 2\frac{r-1}{r^2}$.

Clearly, $F$ is $r$-critical. Also, when removing an edge $xy$ of $F$, we further reduce the chromatic number if and only if $\{x, y\} \cap \{u, v\} = \emptyset$, where $uv$ is the edge removed from $K_{r+2}$. It follows that

$$c(n_1, \ldots, n_r; F) = \sum_{i=2}^{r} \binom{n_i}{2} \prod_{\substack{2 \leq j \leq r \\ j \neq i}} n_j = \frac{(n - n_1 - r + 1)}{2} \prod_{i=2}^{r} n_i.$$

Therefore, $\alpha_F = \frac{r-1}{2r^r}$, $\zeta_F = -\frac{1}{2r^{r-2}}$, and $\pi_F = \frac{r-1}{r^2}$.

On the other hand,

$$P_F(\boldsymbol{\xi}) = \sum_{i=1}^{r} \frac{\xi_i^2}{2} \prod_{\substack{1 \leq j \leq r \\ j \neq i}} \xi_j = \frac{1}{2} \left( \sum_{i=1}^{r} \xi_i \right) \prod_{i=1}^{r} \xi_i.$$

Therefore, if $\sum_{i=1}^{r} \xi_i = \rho + \frac{r-1}{r}$ is fixed, then by convexity $P_F(\boldsymbol{\xi})$ is minimized by picking $\boldsymbol{\xi} = (\rho, 1/r, \ldots, 1/r)$, implying that $\rho_F = \infty$ by Lemma 2.11. Thus Theorem 3.5 implies the desired lower bound on $c_{1,t}(F)$ for every $t$ and it remains to show the upper bound.

Next, we fix an arbitrary integer $t \not\equiv 1 \pmod{r}$ and show that $c_{1,t}(F) \leq \pi_F$. Choose small $\epsilon > 0$ and let all definitions of Section 3.1 apply (with $\delta_0 \ll \epsilon$). Take arbitrary $q \leq (\pi_F + \epsilon)n$. Since $\rho_F = \infty$, the proof of Theorem 3.5 (more specifically, Inequality (18)) shows that, for all sufficiently large $n$, every graph in $\mathcal{H}_F(n, q)$ satisfies $M = \emptyset$. Therefore, we need only to compare graphs obtained from a complete $r$-partite graph by adding extra edges.

In what follows, we identify a graph $H^* \in \mathcal{T}_r(n,q)$ for which $\#F(H^*) = t_F(n,q)$. We then show that we can beat $H^*$ for $q \geq (\pi_F + \delta_0)n$ by using a non-equitable partition. In order to simplify our calculations, we assume that $\lceil n/r \rceil$ is even. Clearly, this happens for infinitely many values $n \equiv t \pmod{r}$, so if we can beat $H^*$ for such $n$, this still implies that $c_{1,t}(F) \leq \pi_F$.

Let $V_1, \ldots, V_r$ be the parts of the Turán graph $T_r(n)$. Assume that $n_1 \geq \cdots \geq n_r$ where $n_i = |V_i|$. Form $H^* \in \mathcal{T}_r(n,q)$ by placing all $q$ extra edges in $V_1$ so that the corresponding bad graph $B^*[V_1]$ is triangle-free and *almost-regular* (that is, the $B^*$-degrees of any two vertices of $V_1$ differ at most by 1). An example of $H^*$ can be obtained by letting $B^*$ form a matching plus isolated vertices when $q \leq n_1/2$ and a path plus disjoint edges otherwise. (Note that $q \leq (\pi_F + \epsilon)n < n_1$ so we have that $\Delta(B^*) \leq 2$.)

**Claim 4.1.** *If $n_1$ is even, then $t_F(n,q) = \#F(H^*)$.*

*Proof of Claim.* For $H \in \mathcal{T}_r(n,q)$ with $E(H) \supseteq E(K(V_1, \ldots, V_r))$, let $\#_2 F(H)$ be the number of $F$-subgraphs that use at most 2 edges of $E(H) \setminus E(K(V_1, \ldots, V_r))$.

Take an arbitrary $H \in \mathcal{T}_r(n,q)$. We need to show that $\#F(H^*) \leq \#F(H)$. Since $B^*$ has no triangle, we have that $\#_2 F(H^*) = \#F(H^*)$. Thus it is enough to show that

$$\#_2 F(H^*) \leq \#_2 F(H),$$

and we can assume that $H$ minimizes $\#_2 F(H)$ among all graphs in $\mathcal{T}_r(n,q)$.

For $i \in [r]$ let $B_i = H[V_i]$. Note that a pair of edges from two different parts $V_i$ and $V_j$ contributes $\left(\binom{r+2}{2} - 2\right) \prod_{h \neq i, j} n_h$ to $\#_2 F(H)$: we have to pick one vertex from every other part, obtaining a copy of $K_{r+2}$, and then remove an edge different from the two initial edges. Also, if a copy of $F$ uses exactly 2 extra edges from the same part $V_i$, then these edges have to be adjacent and this pair contributes $\prod_{j \neq i} n_j$ to $\#_2 F(H)$. Thus if we replace $B_i$ by an arbitrary graph $B'_i$ of the same size, then the obtained graph $H'$ satisfies

$$\#_2 F(H') - \#_2 F(H) = \sum_{x \in V_i} \left(\binom{\deg_{B'_i}(x)}{2} - \binom{\deg_{B_i}(x)}{2}\right) \prod_{\substack{1 \leq j \leq r \\ j \neq i}} n_j.$$

The strict convexity of $\binom{m}{2}$ over integer $m \geq 0$ implies that $B_i$ is almost-regular.

Let $i \in [r]$ be such that $|B_i|$ is maximal. Let us show that we can assume that $i = 1$. If $n_i = n_1$, then we can just swap the labels of the parts $V_1$ and $V_i$. Otherwise, we have $n_i < n_1$. Since $2|B_1| \leq |B_1| + |B_i| \leq q < n_1$, there is $x \in V_1$ which is an isolated vertex of $B_1$. Let $H'$ be obtained from $H$ by 'moving' $x$ from $V_1$ to $V_i$. Formally, define $U_1 = V_1 \setminus \{x\}$, $U_i = V_i \cup \{x\}$, and $U_j = V_j$ for $j \in [r] \setminus \{1, i\}$ and let $H'$ be the (edge-disjoint) union of $K(U_1, \ldots, U_r)$ and $\cup_{j \in [r]} B_j$. If we compare the edge sets of $H'$ and $H$, then we just removed all edges between $x$ and $V_i$ but added all edges between $x$ and $V_1$. Since $n_i = n_1 - 1$, the obtained graph $H'$ is still in $\mathcal{T}_r(n,q)$. Let us consider $\#_2 F(H') - \#_2 F(H)$, where $\#_2 F(H')$ is calculated with respect to the parts $U_j$. By cancellations, it is enough to consider only those $F$-subgraphs that use at least one changed edge at $x$. Since $(|V_1|, |V_i|) = (|U_i|, |U_1|)$, the contributions from $F$-subgraphs that do not use any edge from $B_1 \cup B_i$ also cancel each other. Therefore, $\#_2 F(H') - \#_2 F(H)$ can be determined by looking at those $F$-subgraphs that use $x$ and at least one bad edge from $B_1 \cup B_i$:

$$\#_2 F(H') - \#_2 F(H) = \left(\sum_{u \in V_1} \binom{\deg_{B_1}(u)}{2} - \sum_{u \in V_i} \binom{\deg_{B_i}(u)}{2}\right) \prod_{\substack{2 \leq j \leq r \\ j \neq i}} n_j$$

$$+ \; (|B_1| - |B_i|) \sum_{\substack{2 \leq j \leq r \\ j \neq i}} |B_j| \left(\binom{r+2}{2} - 2\right) \prod_{\substack{2 \leq h \leq r \\ h \neq i, j}} n_h.$$

Since $B_1$ is almost-regular and contains a vertex of degree 0, all its degrees are at most 1. Thus the first summand in the above formula for $\#_2 F(H') - \#_2 F(H)$ is non-positive. By our assumptions, $|B_1| \leq |B_i|$ so the second summand is non-positive too. By the optimality of $H$, we conclude that $\#_2 F(H') = \#_2 F(H)$. Now, redefine $H$ to be $H'$ (and relabel $U_i$ to be the first part). Thus we can indeed assume that $\max_{i \in [r]} |B_i| = |B_1|$.

Suppose that some $B_i$ with $i \ge 2$ is non-empty for otherwise we are done (since, as we have already argued, $B_1$ is almost-regular).

If the graph $B_1$ on $V_1$ has at least two isolated vertices, $x, y \in V_1$, then by removing some $uv \in B_i$ and adding $xy$ to $B_1$, we obtain a graph $H' \in \mathcal{T}_r(n,q)$ with

$$\#_2 F(H') - \#_2 F(H) \le (|B_i| - 1 - |B_1|) \sum_{\substack{2 \le j \le r \\ j \ne i}} |B_j| \left( \binom{r+2}{2} - 2 \right) \prod_{\substack{2 \le h \le r \\ h \ne i,j}} n_h \le 0.$$

(Note that the first inequality holds because $\zeta_F < 0$ and $n_1 \ge n_i$ imply that the number of $F$-subgraphs with $xy$ as the only bad edge is at most that for $uv$.) By iteratively replacing $H$ with $H'$, we can assume that $B_1$ has no isolated vertices. (Recall that $n_1$ is even.) In particular, it follows that $2|B_1| \ge n_1$.

Now, we are ready to compare $H$ directly with $H^*$. In $H$, every edge $uv \in B_j \setminus B_1$ forms at least

$$c(n,F) + |B_1| \left( \binom{r+2}{2} - 2 \right) \prod_{\substack{2 \le h \le r \\ h \ne j}} n_h$$

copies of $F$. Since $B_1$ is almost-regular, it has exactly $2|B_1| - n_1$ vertices of degree 2 while all other vertices have degree 1. As each vertex of degree 2 gives $\prod_{j=2}^{r} n_j$ copies of $F$ that use both edges incident to it, it follows that

$$\#_2 F(H) \ge qc(n,F) + |B_1| \sum_{j=2}^{r} |B_j| \left( \binom{r+2}{2} - 2 \right) \prod_{\substack{2 \le h \le r \\ h \ne j}} n_h + (2|B_1| - n_1) \prod_{j=2}^{r} n_j. \qquad (26)$$

On the other hand, $H^*$ has $qc(n,F)$ copies of $F$ that use exactly one bad edge while there are $2q - n_1$ vertices of $B^*$-degree 2. Thus

$$\#_2 F(H^*) = qc(n,F) + (2q - n_1) \prod_{j=2}^{r} n_j. \qquad (27)$$

By subtracting (27) from (26) and using first that $q = \sum_{i=1}^{r} |B_i|$ and $n_1 = \max_{i \in [r]} n_i$, and then that $\binom{r+2}{2} - 2 \ge 4$ and $2|B_1| \ge n_1$, we obtain that

$$\#_2 F(H) - \#_2 F(H^*) \ge \sum_{j=2}^{r} |B_j| \left( \left( \binom{r+2}{2} - 2 \right) |B_1| - 2n_1 \right) \prod_{\substack{2 \le h \le r \\ h \ne j}} n_h \ge 0.$$

This finishes the proof of Claim 4.1. ∎

Let $q$ satisfy $\delta_0 n \le q - \pi_F n \le \epsilon n$. We now compare the graph $H^*$ of Claim 4.1 with a graph $H$ that is obtained from $K(U_1, \ldots, U_r)$ with part sizes $|U_1| = n_1 + 1$, $|U_2| = n_2 - 1$, and $|U_i| = n_i$ for $i \ge 3$ by placing $q + 1$ extra edges into $U_1$ to form an almost-regular triangle-free graph. Since $n \not\equiv 1 \pmod{r}$, we have that $n_1 = n_2$ and thus $|H| = |H^*| = t_r(n) + q$. Also,

$$\#F(H) \le (q+1)(c(n,F) + \zeta_F n^{r-1}) + (2q - n_1)(n/r)^{r-1} + O(n^{r-1}).$$

Thus, we have that

$$\#F(H) - \#F(H^*) \le \alpha_F n^r + (\pi_F + \delta_0) \zeta_F n^r + o(n^r) < \delta_0 \zeta_F n^r / 2 < 0.$$

This and Claim 4.1 imply the upper bound $c_{1,t}(F) \le \pi_F$ for all $t \not\equiv 1 \pmod{r}$, as required.

It remains to prove that $c_{1,1}(F) \le 2\pi_F$. Let $\ell$ be large and let $n = r\ell + 1$. Let $q = (2\pi_F \pm \epsilon)n$ be given. We first determine $t_F(n,q)$ by constructing a graph $H^* \in \mathcal{T}_r(n,q)$ such that $\#F(H^*) = t_F(n,q)$. Then, for $q \ge (2\pi_F + \delta_0)n$, we exhibit a different graph which beats this bound.

Let $V_1 \cup \ldots \cup V_r$ be the parts of the Turán graph $T_r(n)$ where $|V_1| = \ell + 1$ and $|V_i| = \ell$ for $2 \le i \le r$. As before, we form the graph $H^* \in \mathcal{T}_r(n,q)$ by placing all $q$ extra edges in $V_1$ in such a way that the corresponding bad graph is triangle-free and almost-regular.

**Claim 4.2.** $t_F(n,q) = \#F(H^*)$.

*Proof of Claim.* Let $H \in \mathcal{T}_q(n,r)$. As $H^*$ contains no copy of $F$ using three or more bad edges, we once again consider $\#_2F(H)$, the number of copies of $F$ in $H$ that use at most 2 bad edges, and show that $\#_2F(H^*) \leq \#_2F(H)$. We, therefore, assume that $H$ minimizes $\#_2F(H)$ among all graphs in $\mathcal{T}_q(n,r)$.

For $i \in [r]$, let $B_i = H[V_i]$. As argued before, it follows by convexity that $B_i$ is an almost-regular graph for all $i \in [r]$. Relabel the parts $V_2, \ldots, V_r$ so that $|B_i| \geq |B_j|$ whenever $2 \leq i < j$.

Additionally, we can assume that $|B_1| \geq |B_2|$. In order to show this, suppose that $|B_2| > |B_1|$. Create a new graph $H'$ from $H$ by replacing $H[V_1]$ and $H[V_2]$ by triangle-free almost-regular graphs $B_1'$ and $B_2'$ such that $|B_1'| = |B_2|$ and $|B_2'| = |B_1|$. It is enough to show that $\#_2F(H') \leq \#_2F(H)$. In calculating $\#_2F(H) - \#_2F(H')$, we need only consider the change in copies of $F$ involving edges in $B_1, B_2, B_1'$, and $B_2'$. These come in four flavors:

1. Copies using exactly one bad edge: as $\zeta_F < 0$, there are more such copies in $H$ than in $H'$.
2. Copies using one edge each from parts $V_1$ and $V_2$: $H$ and $H'$ contain an equal number of such copies as $|B_1||B_2| = |B_1'||B_2'|$.
3. Copies using a pair of adjacent bad edges: Let $p_1, p_2, p_1', p_2'$ denote the number of such pairs in $B_1, B_2, B_1', B_2'$, respectively. Their contributions to copies of $F$ in $H$ and $H'$ is then given by

$$(p_1 - p_1') \prod_{i \neq 1} n_i + (p_2 - p_2') \prod_{i \neq 2} n_i = (n_2(p_1 - p_1') + n_1(p_2 - p_2')) \prod_{i \neq 1,2} n_i$$
$$= (n_2(p_1 + p_2 - p_1' - p_2') + (n_1 - n_2)(p_2 - p_2')) \prod_{i \neq 1,2} n_i.$$

This is non-negative because $p_2 \geq p_2'$ (as $|B_2| > |B_2'|$) and $p_1' + p_2' \leq p_1 + p_2$. To see the last inequality, view the transition from $H$ to $H'$ as an iterative process where each step increases $|B_1|$ by 1 and decreases $|B_2|$ by 1. One such step increases $p_1$ by two lowest degrees of $H[V_1]$ and decreases $p_2$ by at least two lowest degrees of $H[V_2]$ after the edge removal. Thus it cannot increase $p_1 + p_2$.

4. Copies using one bad edge from $V_1$ or $V_2$, and one bad edge from a different part: Here, the difference in the number of such copies of $F$ containing a bad edge $e \in B_i$, $i \neq 1,2$, is

$$(|B_1| - |B_1'|)\left(\binom{r+2}{2} - 2\right) \prod_{j \neq 1,i} n_j + (|B_2| - |B_2'|)\left(\binom{r+2}{2} - 2\right) \prod_{j \neq 2,i} n_j$$
$$= \left(\binom{r+2}{2} - 2\right)(|B_2| - |B_1|)(n_1 - n_2) \prod_{j \neq 1,2,i} n_j > 0.$$

As a result, we can indeed assume that $|B_1| \geq |B_2|$.

If $B_2$ (and thus each of $B_i$ for $i \geq 3$) is empty, then there is nothing do to, so suppose otherwise. Furthermore, note that $q = (2\pi_F \pm \epsilon)n = \left(2 - \frac{2}{r} \pm 2r\epsilon\right)\ell$ and, thus, the almost-regular graph $B_1$ contains at least $\ell/2r$ vertices of degree at most 3. We split the rest of the argument depending on the value of $r$.

First, let $r \geq 3$. Create a new graph $H'$ from $H$ by removing a bad edge from $B_2$ and adding a new edge between two vertices of degree at most 3 and at distance larger than 2 in $H[V_1]$. Since $\zeta_F < 0$ and $|V_1| \geq |V_2|$, the number of $F$-subgraphs that use the moved edge as the only bad edge cannot increase. This new edge will form $\ell^{r-1}$ copies of $F$ with each of the at most six adjacent bad edges. However, the number of copies of $F$ that use two bad edges from two different parts decreases by at least $(|B_1| - |B_2| + 1)\left(\binom{r+2}{2} - 2\right)\ell^{r-1} > 6\ell^{r-1}$. This contradicts our choice of $H$ as the minimizer of $\#_2F$ over the set $\mathcal{T}_q(n,r)$.

Let $r = 2$. Here $q = (2\pi_F \pm \epsilon)n = (1 \pm 3\epsilon)\ell$. Also, note that

$$c(\ell, \ell+1; F) = \binom{\ell+1}{2} = \binom{\ell}{2} + \ell = c(\ell+1, \ell; F) + \ell.$$

As $H^*$ has all $q$ edges in $B_1^*$, we have $\#F(H^*) \leq qc(n,F) + \ell^2 + 12\epsilon\ell^2$. Meanwhile,

$$\#_2F(H) \geq qc(n,F) + |B_2|\ell + 2(q - \ell/2 - |B_2|)\ell + 4(q - |B_2|)|B_2|,$$

where the third term counts copies containing a pair of adjacent edges in $B_1$ and the fourth term counts copies using one edge each from $B_1$ and $B_2$. The lower bound on $\#_2 F(H) - \#F(H^*)$ that we obtain is clearly increasing in $q$. Since $q \geq (1-3\epsilon)\ell$, we get that $\#_2 F(H) - \#F(H^*) \geq Q(|B_2|)$, where $Q(x) = -4x^2 + 3\ell(1-4\epsilon)x - 18\epsilon\ell^2$ is a quadratic polynomial that is concave and symmetric around $x_0 = (\frac{3}{8} - \frac{3}{2}\epsilon)\ell$. Routine calculations show that, for example, $Q(7\epsilon\ell) > 0$. Thus $|B_2|$ is either at most $7\epsilon\ell$ or at least $2x_0 - 7\epsilon\ell$. However, the latter is impossible since $|B_2| \leq q/2 \leq (1+3\epsilon)\ell/2$ and $\epsilon$ is small. Thus necessarily $|B_2| < 7\epsilon\ell$. Now, move an edge from $B_2$ to $B_1$, so that the new edge is adjacent to at most 4 bad edges and creates no triangle. The difference in the $F$-count is at most $-\ell - 4(|B_1| - |B_2| + 1) + 4\ell$. The last expression is negative since $|B_1| - |B_2| = q - 2|B_2| > (1-3\epsilon)\ell - 14\epsilon\ell$. This contradicts our assumption that $H$ minimizes $\#_2 F(H)$ over all graphs in $\mathcal{T}_q(n, r)$ and shows that $B_2 = \emptyset$, as required. This finishes the proof of Claim 4.2. ∎

What remains to show now is that $H^*$ can be beaten by perturbing the sizes of the vertex sets. Construct a graph $H$ with vertex set $V = V_1 \cup \ldots \cup V_r$, where $|V_1| = \ell + 2$, $|V_2| = \ell - 1$ and $|V_i| = \ell$ for $3 \leq i \leq r$, and place $q + 2$ edges in $V_1$ such that $B_1 = H[V_1]$ is triangle-free and almost-regular. Note that as $B_1$ has one more vertex and 2 more edges than does $B_1^*$, it has at most 6 more adjacent pairs of bad edges than does $B_1^*$. So, $\#F(H^*) - \#F(H)$ is at least

$$
\begin{aligned}
qc(n, F) &- (q+2)c(\ell+2, \ell-1, \ell, \ldots, \ell; F) - 6(\ell-1)\ell^{r-2} \\
&= q(r-1)\binom{\ell}{2}\ell^{r-2} - (q+2)\left(\binom{\ell-1}{2}\ell^{r-2} + (r-2)\binom{\ell}{2}(\ell-1)\ell^{r-3}\right) - 6(\ell-1)\ell^{r-2} \\
&= \frac{1}{2}(\ell-1)\ell^{r-2}(qr - 2r\ell + 2\ell + 2r - 12).
\end{aligned}
$$

This quantity is positive when $q \geq 2\frac{r-1}{r}\ell + 5$, thus proving that $c_{1,1}(F) \leq 2\pi_F$.

### 4.2. Non-tightness of Theorem 3.5

We now exhibit a graph for which $c_1(F) > \min(\pi_F, \rho_F)$ and some complicated phenomena occur at the threshold, as described at the end of the Introduction. Let $F$ be the graph in Figure 2. As we will see, for this graph, $\rho_F = \hat{\rho}_F = \infty$, so we have to show that $c_1(F) > \pi_F$. Roughly, this inequality is strict because, for not too large $q$, we can reduce the number of copies of $F$ by distributing the bad edges among the two parts of $K_{n/2,n/2}$ instead of placing them all into one part.

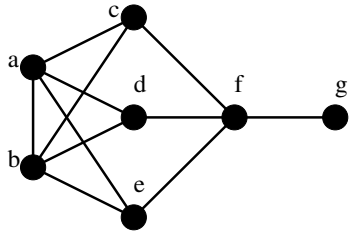

Figure 2: Example for non-tightness of Theorem 3.5.

**Theorem 4.3.** *The graph $F$ of Figure 2 satisfies $c_{1,0}(F) = \frac{3-\sqrt{5}}{4}$ and $c_{1,1}(F) = 1/3$ (while $\pi_F = 1/6$ is strictly less than $c_1(F) = \min(c_{1,0}(F), c_{1,1}(F))$).*

*Proof.* First note that $F$ is 2-critical and $ab$ is the unique critical edge. There is a unique (up to isomorphism) 2-coloring $\chi$ of $F - ab$ with $\chi^{-1}(1) = \{a, b, f\}$ and $\chi^{-1}(2) = \{c, d, e, g\}$. It readily follows that

$$c(n_1, n_2; F) = (n_1 - 2)\binom{n_2}{3}(n_2 - 3)$$

and $\alpha_F = (3! \cdot 2^5)^{-1}$. Taking derivatives, we observe that $\zeta_F = -2^{-5}$ and $\pi_F = 1/6$.

We also have
$$P_F(\boldsymbol{\xi}) = \frac{1}{4 \cdot 3!}(\xi_1 \xi_2^3 + \xi_1^3 \xi_2),$$
which, if we fix $\xi_1 + \xi_2$, is minimized by maximizing the difference. It follows (for example, from Lemma 2.11) that $\rho_F = \infty$.

Note that there exists a 2-coloring $\chi^*$ of $F - ab - fg$ with $\chi^*(a) = \chi^*(b) = \chi^*(f) = \chi^*(g) = 1$. Furthermore, if $u_1 v_1, u_2 v_2$ are two distinct edges in $F$, there is no 2-coloring $\chi''$ of $F - u_1 v_1 - u_2 v_2$ with $\chi''(u_1) = \chi''(v_1)$ and $\chi''(u_2) = \chi''(v_2)$ unless $\{u_1 v_1, u_2 v_2\} = \{ab, fg\}$ and $\chi''$ is isomorphic to $\chi^*$. That is, for any $H \in \mathcal{T}_r(n,q)$, the only copies of $F$ in $H$ that use exactly two bad edges correspond to $\chi^*$.

Take any $H \in \mathcal{H}_F(n,q)$ with $q = O(n)$. Fix arbitrary $\epsilon > 0$. Let the definitions of Section 3.1 apply. Once again, as $\rho_F = \infty$, we conclude from (18) that $M = \emptyset$, that is, there are no missing edges. Assume without loss of generality that $|V_1| \geq |V_2|$. By Claim 3.6, we have
$$|B_1| \geq |B_1| - |B_2| \geq (1 - \delta_0)(|V_1| - |V_2| - 1)\pi_F n.$$

First, we prove the claimed lower bound on $c_{1,0}(F)$. Suppose that $n$ is even and $q < (3 - \sqrt{5} - \epsilon)n/4$. We derive a contradiction by assuming that $|V_1| \geq |V_2| + 2$ (and that $\delta_0 \ll \epsilon$).

Let us show that $B_2 = \emptyset$. If $B_2 \neq \emptyset$, an edge $uv \in B_2$ is contained in at least $c(n_2, n_1; F) > c(n_1, n_2; F) - 2\zeta_F n^4 - \delta_0 n^4$ copies of $F$. However, if we remove $uv$ and replace it with an edge $xy$ where $x, y \in V_1$ have $d_B(x), d_B(y) \leq O(1)$, we form at most $c(n_1, n_2; F) + q \times 4\binom{n/2}{3} + \delta_0 n^4$ copies of $F$. As
$$-2\zeta_F n^4 = 2^{-4} n^4 > q \times 4\binom{n/2}{3} + \delta_0 n^4,$$
this alteration reduces the number of copies of $F$. Since $H$ is optimal, we conclude that $B_2 = \emptyset$, as claimed.

In addition, the maximum degree in $B_1$ is at most $\delta_0 n$ (since a vertex of $B_1$-degree at least $\delta_0 n$ creates at least $\delta_1 n^6$ copies of $F$ that use more than one bad edge while, if $\Delta(B_1) = O(1)$, then we would have only $O(n^5)$ such copies).

We conclude that $B_1$ contains at least $(1 - \delta_0)q^2/2$ disjoint pairs of edges, each of which forms $4\binom{|V_2|}{3}$ copies of $F$. It follows that
$$\#F(H) \geq (q + a^2)(c(n, F) + a\zeta_F n^4) + q^2 \frac{n^3}{24} - \delta_0 n^5, \tag{28}$$
where $a = |V_1| - n/2 = n/2 - |V_2| \geq 1$. Since $a = o(n)$, the main terms involving $a$ are $a^2 n^5/(3!\, 2^5) - aqn^4/2^5$. Since e.g. $q < n/5$, the last expression is minimized when $a = 1$ with the minimum being by $\Omega(n^5)$ smaller than any other value for $a \geq 2$. Thus (28) is minimized when $a = 1$.

On the other hand, construct $H^* \in \mathcal{T}_2(n,q)$, where we place $q/2$ edges in each of $B_1$ and $B_2$, thereby forming at most $q^2/4$ pairs of bad edges that lie in the same part. Thus,
$$\#F(H^*) \leq qc(n, F) + \frac{q^2}{4} \times 4\binom{n/2}{3} + \delta_0 n^5. \tag{29}$$

Comparing the above quantities, we get that $\#F(H) > \#F(H^*)$, a contradiction that proves the stated lower bound on $c_{1,0}(F)$. (Note that the right-hand sides of (28), for $a = 1$, and (29) become equal when $q$ is around $\frac{3-\sqrt{5}}{4}n$.)

The upper bound on $c_{1,0}(F)$ follows by noting that inequalities (28) and (29) may be 'reversed' by replacing the last term with $\pm \delta_0 n^5$, where $H^* \in \mathcal{H}_F(n,q)$ is arbitrary and $H$ is obtained by adding $q$ extra edges into the larger part of $K(V_1, V_2)$ with $|V_1| = |V_2| + 2$.

Finally, let us briefly discuss the case when $n = 2\ell + 1$ is odd (and $q \leq (1/3 + \epsilon)n$). As before, we can assume that $M = \emptyset$. Thus we have to identify an asymptotically optimal way of adding edges to $K(V_1, V_2)$, where $|V_1| = \ell + 1 + a$ and compare the cases $a = 0$ and $a \geq 1$. Let $q_i = |B_i|$ be the number of edges inside $V_i$ for $i = 1, 2$.

First, let $a = 0$. We have $q_1 + q_2 = q$. It is advantageous to spread the bad edges inside each $V_i$ uniformly; then the number of $F$-subgraphs is

$$(q_1 + q_2)\alpha_F n^5 - q_2 \zeta_F n^4 + 4\left(\binom{q_1}{2} + \binom{q_2}{2}\right)\binom{n/2}{3} \pm \delta_0 n^5.$$

For $q = q_1 + q_2 < 3n/8 - \epsilon$, the main terms are minimized by letting $q_2 = 0$, which gives $\alpha_F q n^5 + q^2 n^3/24$.

If $a \geq 1$, then $q_1 + q_2 = q + a^2 + a$. It is routine to see that the optimal way is to let $a = 1$ and $q_2 = 0$, giving $(q+2)(\alpha_F n^5 + \zeta_F n^4) + 4\binom{q}{2}\binom{n/2}{3} \pm \delta_0 n^4$ copies of $F$. After removing the main term $q\alpha_F n^5$ from both expressions, we have to compare terms of order $n^5$, that is, $q^2 n^3/24$ versus $n^5/(3 \cdot 2^5) - q n^4/2^5 + q^2 n^3/12$. One can see that they become equal when $q = n/3$ and conclude that $c_{1,1}(F) = 1/3$, completing the proof. □

### 4.3. Pair-free graphs

One property of the graph in Figure 2 is that there exists a 2-coloring of the vertices that would be a proper 2-coloring with the deletion of exactly two edges from one color class. We now consider graphs which do not have this property.

**Definition 4.4.** Let $F$ be an $r$-critical graph. We say that $F$ is *pair-free* if there do not exist two (different, but not necessarily disjoint) edges $u_1 v_1, u_2 v_2$ and a proper $r$-coloring $\chi$ of $F - u_1 v_1 - u_2 v_2$ such that $\chi(u_1) = \chi(u_2) = \chi(v_1) = \chi(v_2)$.

Many interesting graphs belong to this class, e.g., odd cycles and cliques. In addition, graphs obtained from the complete $r$-partite graph $K_{s_1,\ldots,s_r}$ by adding an edge to the part of size $s_1$ are pair-free if $s_i \geq 3$ for all $i \geq 2$.

**Proposition 4.5.** Let $F$ be pair-free and let $t = -\operatorname{sign}(\zeta_F)$. Then $c_{1,t}(F) \leq 2\pi_F$ and $c_{1,i}(F) \leq \pi_F$ for $i \not\equiv t \pmod{r}$.

*Proof.* We prove the case $n \not\equiv t \pmod{r}$; the other case follows in a similar manner. Let $n$ be large and $q = (\pi_F + \delta_0)n$. Write $n = n_1 + \ldots + n_r$, where $c(n, F) = c(n_1, \ldots, n_r; F)$ and the sequence $(n_1, \ldots, n_r)$ is monotone. Since $n \not\equiv t \pmod{r}$, we have $n_1 = n_2$. Consider the partition $n = n'_1 + n'_2 + \ldots + n'_r$ where $n'_1 = n_1 + t$, $n'_2 = n_2 - t$ and $n'_i = n_i$ for $i = 3, \ldots r$. Construct $H'$ as follows: start with $K(V'_1, \ldots, V'_r)$, where $|V'_i| = n'_i$, and place $q + 1$ edges in $V'_1$ to form an almost regular bipartite graph. We claim that $\#F(H') < \#F(H)$ for any $H \in \mathcal{T}_r(n, q)$.

In $H'$, each bad edge is contained in at most $c(n, F) - |\zeta_F| n^{f-3} + O(n^{f-4})$ copies of $F$ that contain no other bad edge. As $F$ is pair-free, no copy of $F$ contains exactly two bad edges. In addition, we may bound the number of copies of $F$ that use at least three bad edges by $O(n^{f-3})$.

On the other hand, $\#F(H) \geq qc(n, F)$. Therefore,

$$\begin{aligned}
\#F(H') - \#F(H) &\leq (q+1)\bigl(c(n,F) - |\zeta_F|n^{f-3}\bigr) + O(n^{f-3}) - qc(n,F) \\
&< \alpha_F n^{f-2} - (\pi_F + \delta_0)n|\zeta_F|n^{f-3} + O(n^{f-3}) < 0,
\end{aligned}$$

proving the proposition. □

For odd cycles, this implies that $c_1(C_{2k+1}) = c_{1,0}(C_{2k+1}) = 1/2$ by Theorem 3.5. In fact, with more effort, it is possible to show that $c_{1,1}(C_{2k+1}) = 1$. However, the proof is quite involved as one has to account for copies of $C_{2k+1}$ that may appear in various configurations. We direct the interested reader to [25].

**Corollary 4.6.** If an $r$-critical pair-free graph $F$ satisfies $\deg(P_F) \geq r + 1$ and $\rho_F = \hat{\rho}_F$, then $c_{1,t}(F) = \min(2\pi_F, \hat{\rho}_F)$ and $c_{1,i}(F) = \min(\pi_F, \hat{\rho}_F)$ for every residue $i$ different from $t$, where $t = -\operatorname{sign}(\zeta_F)$.

*Proof.* Since $\deg(P_F) \geq r + 1$, Theorem 3.5 gives the stated lower bounds. On the other hand, the proof of Theorem 3.12 does not rely on the residue of $n$ modulo $r$ and shows that $c_{1,i}(F) \leq \hat{\rho}_F$ for every $i$. Combining this with Proposition 4.5 we obtain the desired upper bound. □

Our next result shows that the assumptions of Corollary 4.6 hold for $F = K^+_{s,t}$ where $t \geq 3$.

**Lemma 4.7.** *Let $s, t \geq 2$ and $F = K_{s,t}^+$ be obtained from the complete bipartite graph $K_{s,t}$ by adding an edge to the part of size $s$. Then $\rho_F = \hat{\rho}_F$. Furthermore, if $t = 2, 3$, then $\rho_F = \infty$.*

*Proof.* Assume that $\rho_F < \infty$ is finite as otherwise $\rho_F = \hat{\rho}_F = \infty$ and we are done.

Clearly, $F = K_{s,t}^+$ is 2-critical and
$$c(n_1, n_2; K_{s,t}^+) = \binom{n_2}{t}\binom{n_1-2}{s-2}.$$

It readily follows that
$$\alpha_F = \frac{2^{-(t+s-2)}}{t!(s-2)!}, \quad \zeta_F = (s-t-2)\frac{2^{-(t+s-3)}}{t!(s-2)!}, \quad \text{and}$$

$$\pi_F = \begin{cases} \infty, & \text{if } t = s-2, \\ (2\,|t-s+2|)^{-1}, & \text{otherwise.} \end{cases}$$

On the other hand,
$$P_F(\boldsymbol{\xi}) = \frac{2^{-s+2}}{t!(s-2)!}(\xi_1 \xi_2^t + \xi_1^t \xi_2).$$

As $P_F$ is a homogeneous polynomial, we restrict ourselves to $\xi_1 + \xi_2 = 1$. Namely, let
$$\varphi_{s,t}(y) = P_F(1/2 + y, 1/2 - y) = \frac{2^{-s+2}}{t!(s-2)!}\Big((1/2+y)(1/2-y)^t + (1/2+y)^t(1/2-y)\Big).$$

We observe that $\varphi_{s,t}(y)$ is an even function with $\varphi_{s,t}(1/2) = \varphi_{s,t}(-1/2) = 0$ and $\varphi_{s,t}(0) = \alpha_F$. Routine calculations show that the coefficient $s_k$ of $y^k$ in the derivative $\varphi'_{s,t}(y)$ for $0 \leq k \leq t$ is
$$\Big((-1)^k - 1\Big)(2k + 1 - t)\binom{t}{k}\frac{2^{-s+2+(k-t)}}{t!(s-2)!}.$$

It follows that $s_k = 0$ when $k = (t-1)/2$ or $k$ is even. Otherwise, if $k < (t-1)/2$ (resp. $k > (t-1)/2$), then $s_k$ is positive (resp., negative).

Let us first do the case $t \geq 4$. Since the coefficients of $\varphi'_{s,t}(y)$ change sign exactly once, $\varphi'_{s,t}(y)$ has at most one positive root by Descartes' rule of signs. As $\varphi''_{s,t}(0) = \frac{2^{-s+2}}{t!(s-2)!}\frac{4t(t-3)}{2^t} > 0$ for $t \geq 4$ while the coefficient at the highest degree term of $\varphi'_{s,t}$ is negative, the polynomial $\varphi'_{s,t}(y)$ has exactly one positive root (and, by symmetry, exactly one negative root).

It follows that $(0, \alpha_F)$ is the unique local minimum for $\phi_{s,t}$ with the two roots of $\varphi'_{s,t}$ providing local maxima. Thus the value of $p(\rho)$, as defined by (9), is attained at $\xi_1 = \xi_2 = \frac{\rho}{2} + \frac{1}{4}$ or at $\{\xi_1, \xi_2\} = \{\rho, \frac{1}{2}\}$. The latter solution cannot give $\rho_F$ by Lemma 2.11.

Thus $\rho_F$ is the smallest positive root of the equation $P_F(\frac{\rho}{2} + \frac{1}{4}, \frac{\rho}{2} + \frac{1}{4}) = \alpha_F \rho$. After simplifications, the equation becomes $g(\rho) = 0$, where $g(\rho) = (\rho + \frac{1}{2})^{t+1} - \rho$. Let $\rho_0 = 5^{-1/4} - 1/2$. Note that $g(0) > 0$ while $g(\rho_0) \leq 5^{-5/4} - \rho_0 < 0$. Thus $\rho_F < \rho_0$.

Suppose that $\rho_F \neq \hat{\rho}_F$. Then the two sides of the original equation not only coincide but also are tangent at $\rho_F$, which means that $g'(\rho_F) = 0$. Solving the obtained equation $(t+1)(\rho + \frac{1}{2})^t = 1$, we obtain
$$\rho_F = (t+1)^{-1/t} - 1/2 \geq 5^{-1/4} - 1/2 = \rho_0 > \rho_F,$$

which is the desired contradiction.

Finally, let us consider cases when $t = 2$ or $3$. We have that $\phi_{s,2}(y) = \alpha_F(1 - 4y^2)$ and $\phi_{s,3}(y) = \alpha_F(1 - 16y^4)$. It easily follows that the minimum of $P_F$ over $\mathcal{S}_\rho$ is attained for $\{\xi_1, \xi_2\} = \{\rho, \frac{1}{2}\}$. As before, this cannot give $\rho_F$ by Lemma 2.11 and contradicts our assumption that $\rho_F < \infty$. □

**Corollary 4.8.** *Let $s \geq 2$, $t \geq 3$ and $F = K_{s,t}^+$. Then $c_1(F) = c_{1,0}(F) = \min(\pi_F, \rho_F)$ and $c_{1,1}(F) = \min(2\pi_F, \rho_F)$.*

*Proof.* Note that the 2-critical graph $F$ satisfies $\deg(P_F) = t + 1 > 3$. The lower bound follows from Theorem 3.5. On the other hand, note that $F$ is pair-free. Lemma 4.7 gives that $\rho_F = \hat{\rho}_F$ while we have by the proof of Theorem 3.12 that $c_{1,i}(F) \leq \hat{\rho}_F$ for both $i = 0, 1$. Finally, Proposition 4.5 completes the proof. □

When $t = 2$, the graph $K_{s,t}^+$ is not pair-free and Proposition 4.5 does not apply. We consider these graphs in the next section.

### 4.4. $K_{s,2}^+$

**Theorem 4.9.** *Let $s \geq 2$ and $F = K_{s,2}^+$. Then $c_1(F) = c_{1,0}(F) = \begin{cases} \pi_F, & \text{if } s = 2 \text{ or } s \geq 6, \\ \infty, & \text{if } s \in \{3, 4, 5\}. \end{cases}$*

*In addition, $c_{1,1}(F) = \begin{cases} 2\pi_F, & \text{if } s = 2 \text{ or } s \geq 8, \\ 3/8, & \text{if } s = 7, \\ \infty, & \text{if } s \in \{3, 4, 5, 6\}. \end{cases}$*

*Proof.* By Lemma 4.7, we have that $\rho_F = \infty$. Also, recall that $\alpha_F = (2^{s+1}(s-2)!)^{-1}$, $\zeta_F = (s-4)/(2^s(s-2)!)$, and $\pi_F = (2\,|s-4|)^{-1}$.

The case of $K_{2,2}^+ = K_4 - e$ was covered in Section 4.1. Also, if $s = 4$, then $\rho_F = \pi_F = \infty$ and Theorem 3.5 gives that $c_{1,1}(F) = c_{1,0}(F) = \infty$, as required. It thus remains to consider the cases $s = 3$ or $s \geq 5$. Let all definitions of Section 3.1 apply.

Since $\rho_F = \infty$, the proof of Theorem 3.5 gives that $M(H) = \emptyset$ for any $H \in \mathcal{H}_F(n, q)$ and appropriate $q = q(n)$. Therefore, we need only to compare graphs obtained from a complete 2-partite graph by adding extra edges. As observed earlier, $F = K_{s,2}^+$ is not pair-free. However, any copy of $F$ which contains more than one bad edge must have at least two bad edges that are adjacent to each other. Therefore, all such copies can be avoided if the 'extra' edges form a matching.

We use this observation to first show that $c_{1,0}(F) \leq 1/(2(s-4))$ for $s \geq 6$. Let $n = 2\ell$ be large and let $2(q+1) \leq \ell - 1$. For any $H \in \mathcal{T}_2(n, q)$, we have the lower bound $\#F(H) \geq qc(n, F) = q\binom{\ell}{2}\binom{\ell-2}{s-2}$. We construct the graph $H'$ by taking $K(V_1, V_2)$ with $|V_1| = \ell - 1$ and $|V_2| = \ell + 1$ and placing $q + 1$ edges in $V_1$ so that they form a matching. As all copies of $F$ in $H'$ use exactly one bad edge, we have $\#F(H') = (q+1)\binom{\ell+1}{2}\binom{\ell-3}{s-2}$. It is routinely seen that $\#F(H) > \#F(H')$ for $q \geq \ell/(s-4) = \pi_F n$. (For $s = 6$, this can be improved to e.g. $q \geq \ell/2 - 3$, giving a non-empty intersection with the restriction $2(q+1) \leq \ell - 1$.) Thus indeed $c_{1,0}(F) \leq \pi_F$ for $s \geq 6$. On the other hand, the corresponding lower bound is given by Theorem 3.5.

A similar argument shows that $c_{1,1}(F) \leq 1/(s-4)$ for $s \geq 8$. Namely, given a large odd integer $n = 2\ell + 1$, let $H'$ be constructed by taking $K(V_1, V_2)$ with $|V_1| = \ell - 1$, $|V_2| = \ell + 2$ and placing $q + 2 \leq (\ell - 1)/2$ edges in $V_1$ so that the added edges form a matching. We then have $\#F(H') = (q+2)\binom{\ell+2}{2}\binom{\ell-3}{s-2}$ which is seen to be smaller than $qc(n, F) = q\binom{\ell+1}{2}\binom{\ell-2}{s-2}$ for $q \geq 2\ell/(s-4)$ (for $s = 8$, this can be improved to $q \geq \ell/2 - 4$). Thus indeed $c_{1,1}(F) \leq 2\pi_F$ for $s \geq 8$. The converse inequality follows from Theorem 3.5, as desired.

In the remaining cases, we have to take into account copies of $F$ containing more than one bad edge so the calculations are more complicated.

First, let us consider the case where $n = 2\ell$ is even (and tends to infinity). Let $H$ be some graph obtained from $K(V_1, V_2)$ with $|V_1| = \ell - \operatorname{sign}(\zeta_F)a$ and $|V_2| = \ell + \operatorname{sign}(\zeta_F)a$, by placing $b_1$ edges in $V_1$ and $b_2$ edges in $V_2$. Recall that if a copy of $F$ contains exactly two bad edges, then they are adjacent. Furthermore, two adjacent bad edges, say inside $V_1$, create $|V_2|\binom{|V_2|-1}{s-2}$ copies of $F$. Thus

$$\#F(H) \geq (b_1 + b_2)c(\ell, \ell; F) + \left(1 + O\left(\frac{1+a^2}{\ell}\right)\right)\left((b_2 - b_1)a|\zeta_F|(2\ell)^{s-1} + \tau(H)\frac{\ell^{s-1}}{(s-2)!}\right), \quad (30)$$

where we define $\tau(H) = \sum_{v \in V(H)} \binom{d_B(v)}{2}$ to be the number pairs of adjacent edges lying in some part $V_i$. Furthermore, if the bad graph $B$ has maximum degree $O(1)$, then (30) is in fact equality within an additive error term $O((1+a^2)\ell^{s-1})$.

Returning to the claimed lower bound $c_{1,0}(F) \geq \infty$ for $s \in \{3, 5\}$, suppose on the contrary that there is an $h_F(n,q)$-optimal graph $H$ as above with $a \neq 0$ and $q = b_1 + b_2 - a^2 = O(n)$. By symmetry, assume that $a \geq 1$. Since $h_F(n,q) \leq qc(\ell,\ell;F) + O(n^s)$ and $b_1 + b_2 = q + a^2$, we necessarily have that $a = o(\sqrt{n})$.

Next, we reorganize the bad edges in such a way as to minimize $\tau(H)$ while keeping the number of bad edges in each part unchanged. Namely, we construct such a graph (call it $H'$) by letting, for $i = 1, 2$, the corresponding bad graph $B'_i$ have all degrees $d_i$ except for $r_i$ vertices of degree $d_i + 1$, where the integers $d_i$ and $r_i$ come from the integer division of $2b_i$ by $|V_i|$, namely, $2b_i = d_i|V_i| + r_i$ with $0 \leq r_i < |V_i|$. We have by convexity that $\tau(H') \leq \tau(H)$ and, since (30) is equality for $H'$, that $\#F(H') \leq \#F(H) + O(a^2 \ell^{s-1})$. Note that by Claim 3.6

$$b_1 - b_2 \geq (1 - \delta_2)(2a - 1)\pi_F n = 2(1 - \delta_2)(2a - 1)\pi_F \ell. \tag{31}$$

We consider modifying $H'$ by moving $\delta_1 \ell$ edges from $B'_1$ to $B'_2$ (and updating the new bad graphs to be almost-regular). Let us show that the above transformation decreases the value of

$$\lambda a(|B'_2| - |B'_1|) + \tau(H'), \tag{32}$$

by at least $4\delta_1(b_1 - b_2) - \delta_1(2\lambda a + 2 + \delta_0)\ell$, where $\lambda = 2^{s-1}|\zeta_F|(s-2)!$. First, suppose that $r_1 \geq 2\delta_1 \ell$ and $r_2 \leq |V_2| - 2\delta_1 \ell$. Here, $2\delta_1 \ell$ degrees in $B'_1$ decrease from $d_1 + 1$ to $d_1$ and $2\delta_1 \ell$ $B'_2$-degrees increase from $d_2$ to $d_2 + 1$. Since $d_1$ and $d_2$ are bounded by e.g. $5q/n = O(1)$, the value of (32) decreases by

$$-2\delta_1 \ell \lambda a + 2\delta_1 \ell(d_1 - d_2) + O(a^2) = \delta_1(4b_1 - 4b_2 - 2r_1 + 2r_2 - 2\lambda a\ell) + O(a^2). \tag{33}$$

This implies the desired bound since $r_1 < |V_1|$ and $r_2 \geq 0$. If $r_2 > |V_2| - 2\delta_1 \ell$, then we do not have enough vertices of degree $d_2$ in $B'_2$ so the decrease in $\tau(H')$ that we can guarantee is only $\delta_1 \ell(2d_1 - 2(d_2 + 1))$, which is lower by $2\delta_1 \ell$ than the bound used in (33). However, this is essentially compensated by the term $2r_2\delta_1$ in (33) as $r_2 \geq (1 - 3\delta_1)\ell$ now. Likewise, if $r_1 < 2\delta_1 \ell$, then we have to decrease some $B'_1$-degrees from $d_1$ to $d_1 - 1$ and may lose up to $2\delta_1 \ell$ when considering $\tau(H')$. However, we gain at least $2\delta_1(|V_1| - 2\delta_1 \ell)$ when we use the assumption $r_1 < 2\delta_1 \ell$ instead of the previous trivial estimate $r_1 < |V_1|$. We see that the decrease in (32) is as claimed in all cases.

Since $H$ is optimal and (30) is equality if applied to $H'$, the value in (32) cannot decrease by more than $O(a^2)$. Thus

$$b_1 - b_2 \leq \frac{(2^s |\zeta_F|(s-2)!\, a + 2 + 2\delta_0)\ell}{4}. \tag{34}$$

In particular, we immediately see that (31) and (34) cannot be simultaneously satisfied for any $a \geq 1$ when $s = 3$ (then $\zeta_F = -1/8$ and $\pi_F = 1/2$) nor when $s = 5$ (then $\zeta_F = 1/(3 \cdot 2^6)$ and $\pi_F = 1/2$). This contradiction shows that $c_{1,0}(K_{3,2}^+) = c_{1,0}(K_{5,2}^+) = \infty$, as desired.

Let us consider the case when $n = 2\ell + 1$ is odd (and $s \in \{3, 5, 6, 7\}$). Let $n_1$ and $n_2$ satisfy $\{n_1, n_2\} = \{\ell, \ell + 1\}$ and $c(n, F) = c(n_1, n_2; F)$. Thus, $n_1 < n_2$ if and only if $\zeta_F > 0$. Let $H$ be some graph obtained from $K(V_1, V_2)$ with $|V_1| = n_1 - \text{sign}(\zeta_F)a$ and $|V_2| = n_2 + \text{sign}(\zeta_F)a$, by placing $b_1$ edges in $V_1$ and $b_2$ edges in $V_2$. The analog (30) for odd $n$ is

$$\begin{aligned}
\#F(H) &\geq (b_1 + b_2)c(n_1, n_2; F) \\
&+ \left(1 + O\left(\frac{1 + a^2}{\ell}\right)\right)\left(((a+1)b_2 - ab_1)|\zeta_F|(2\ell)^{s-1} + \tau(H)\frac{\ell^{s-1}}{(s-2)!}\right),
\end{aligned} \tag{35}$$

where $\tau(H)$, as before, is the number of pairs of adjacent bad edges. Again, if $\Delta(B) = O(1)$, then this is in fact equality within an additive error term $O((1+a^2)\ell^{s-1})$.

Suppose that we have an optimal $H$ as above that contradicts the theorem. Here part sizes $|V_1|$ and $|V_2|$ differ by at least 3 and $q = b_1 + b_2 - a^2 - a = O(n)$. Because of the symmetry between $a \geq 1$ and $a \leq -2$, let us assume that $a \geq 1$. Again, the main terms show that necessarily $a = o(\sqrt{n})$. Claim 3.6 implies that

$$b_1 - b_2 \geq (1 - \delta_2) \cdot 2a \cdot \pi_F n = 4(1 - \delta_2)a\pi_F \ell. \tag{36}$$

The obvious modification of the argument that led to (34) gives that

$$b_1 - b_2 \leq \frac{(2^{s-1}\,|\zeta_F|\,(s-2)!\,(2a+1) + 2 + \delta_0)\ell}{4}. \tag{37}$$

If $s \in \{3,5\}$, then the obtained bounds (36) and (37) contradict each other, proving that $c_{1,1}(F) = \infty$ then.

Next, suppose that $s = 6$ and $n = 2\ell + 1$ is odd. Here, (36) and (37) do not give a contradiction right away. Namely, they state that

$$(1 - \delta_2)a\ell \leq b_1 - b_2 \leq (2a + 3 + \delta_0)\ell/4.$$

We can only conclude that $a = 1$ so we have to provide an extra argument. In order to get a contradiction for $a = 1$ it is enough to exhibit, for $q = b_1 + b_2 - 2$, some $H^* \in \mathcal{T}_2(n,q)$ with $\#F(H) - \#F(H^*) = \Omega(n^6)$. Let the bad edges of $H^*$ be evenly split among the two vertex classes and form an almost regular graph in each. Also, we can assume that the bad graphs $B_1$ and $B_2$ of $H$ are almost regular. Since $b_1 - b_2 \geq (1-\delta_2)\ell$, if we compare the ordered degree sequences in $B_1$ and $B_2$, then each term in the former is at least by 2 larger than the corresponding term in the latter, except for at most $3\delta_2\ell$ terms. It follows that $\tau(H^*) \leq \tau(H) - (1-\delta_0)\ell$. On the other hand, we know that $b_1 - b_2 \leq (5+\delta_0)\ell/4$. If $b_1 + b_2 = q + 2$ is fixed, then $2b_2 - b_1$ is minimized when $b_1 = (q+2)/2 + (5+\delta_0)\ell/8$ is as large as possible given the above constraints. Thus by (35) we have

$$\#F(H) \geq (q+2)c(n,F) + \left(2 \cdot \frac{4q - (5+\delta_0)\ell}{8} - \frac{4q + (5+\delta_0)\ell}{8} + \tau(H)\right) \frac{\ell^5}{24} + O(\ell^5).$$

Comparing this with $\#F(H^*) = qc(n,F) + (q/2)\ell^5/24 + \tau(H^*)\ell^5/24 + O(\ell^5)$, we obtain

$$\#F(H) - \#F(H^*) \geq 2c(n,F) - \frac{7\ell^6}{192} - \delta_0 \ell^6,$$

which is strictly positive as $c(n,F) = (1/48 + o(1))\ell^6$. This contradicts the optimality of $H$. Therefore, $\mathcal{H}(n,q) \subseteq \mathcal{T}_2(n,q)$ for $q = O(n)$, proving that $c_{1,1}(K_{6,2}^+) = \infty$.

Finally, let show that, for $s = 7$, we have $c_{1,1}(F) = 3/8$. Note that this is strictly larger than $2\pi_F = 1/3$, so we need to prove both bounds. Unfortunately, (36) and (37) do not directly contradict each other here, so some further analysis is needed.

The proof of Theorem 3.5 shows that it is enough to consider graphs without missing edges. Furthermore, a hypothetical counterexample must have at least $2(1-\delta_2)\ell/3$ bad edges by (36). Thus it is enough to restrict ourselves to $q$ with

$$2(1-\delta_2)\ell/3 \leq q \leq (3/4 + \delta_0)\ell. \tag{38}$$

Now consider a graph $H_a$ obtained from $K(V_1, V_2)$ with $|V_1| = \ell - a$ and $|V_2| = \ell + 1 + a$ by adding almost regular graphs $B_1$ and $B_2$ of sizes $b_1 + b_2 = q + a^2 + a$, where $a = o(\sqrt{\ell})$. Given $q$, let us determine the optimal values of $b_1$ and $b_2$. We have

$$\#F(H_a) = (q + a^2 + a)c(n,F) + ((a+1)b_2 - ab_1)\ell^6/80 + \tau(H_a)\ell^6/120 + O(a^2\ell^6).$$

First, if $a = 0$, then the main terms are minimized when $b_1 = \lceil \ell/2 \rceil$, the moment when the addition of an extra edge to $V_1$ starts increasing $\tau(H_0)$ by at least 2. Then $B_2$ is a matching by (38), $\tau(H_0) = O(1)$, and

$$\#F(H_0) = qc(n,F) + (q - \ell/2)\ell^6/80 + O(\ell^6). \tag{39}$$

Next let $a \geq 1$. If $b_1 \leq (\ell - a - 1)/2$, then we can improve the main terms for $\#F(H_a)$ when moving one edge from $B_2$ to $B_1$. So we can assume that $b_1 \geq (\ell - a)/2$. Then, by (38), $B_1$ has only degrees 1 and 2 while $B_2$ is a matching. Thus $\tau(H_a) = 2b_1 - \ell + a = 2b_1 - \ell + O(\sqrt{\ell})$ and

$$\#F(H_a) = (q + a^2 + a)c(n,F) + ((a+1)(q - b_1) - ab_1)\ell^6/80 + (2b_1 - \ell)\ell^6/120 + o(\ell^7).$$

This is smallest when $b_1 = q + a^2 + a$ is maximum possible. By plugging this value of $b_1$, subtracting the estimate (39), and using that $c(n, F) = \ell^7/240 + O(\ell^6)$, we get

$$\#F(H_a) - \#F(H_0) = \frac{((2a^2 + 2a - 1)\ell + (2 - 6a)q - 12a^3)\ell^6}{480} + o(\ell^7).$$

It is routine to see that, if $a \geq 2$, this difference is $\Omega(\ell^7)$ (i.e. positive) for all $q$ as in the allowed interval (38). On the other hand, $\#F(H_1) - \#F(H_0)$ stays $\Omega(\ell^7)$ for $q < (3/4 - \delta_0)\ell$ and becomes negative (of order $\ell^7$) before $q$ reaches $(3/4 + \delta_0)\ell$. Thus indeed $c_{1,1}(K_{7,2}^+) = 3/8$.

This completes the proof of Theorem 4.9. $\square$

**Glossary**

$\alpha_F$ : the leading coefficient (coefficient of $n^{f-2}$) of $c(n,F)$. Page 4

$c_1(F)$ : the threshold constant for when graphs in $\mathcal{T}_r(n,q)$ are optimal. Page 3

$c_{1,i}(F)$ : the threshold constant for when graphs in $\mathcal{T}_r(n,q)$ are optimal where $n \equiv i \pmod{r}$. Page 3

$c_2(F)$ : the threshold constant for when graphs in $\mathcal{T}_r(n,q)$ are asymptotically optimal. Page 2

$c(n,F)$ : the minimum number of copies of $F$ in a graph obtained by adding one edge to $T_r(n)$. Page 2

$c(n_1,\ldots,n_r;F)$ : the number of copies of $F$ obtained by adding one edge to the part of size $n_1$ in the complete $r$-partite graph with parts of size $n_1,\ldots,n_r$. Page 4

$\mathrm{ex}(n,F)$ : The Turán function: maximum number of edges in an $F$-free graph with $n$ vertices. Page 1

$F$ : a graph, typically an $r$-critical graph on $f$ vertices. Page 2

$\mathcal{H}_F(n,q)$ : the set of graphs on $n$ vertices and $\mathrm{ex}(n,F)+q$ edges which contain the smallest number of copies of $F$; the set of optimizers of $h_F(n,q)$. Page 2

$h_F(n,q)$ : the minimum number of $F$-subgraphs in a graph with $n$ vertices and $\mathrm{ex}(n,F)+q$ edges. Page 1

$P_F(\boldsymbol{\xi}) := \frac{1}{\mathrm{Aut}(F)} \sum_{u \text{ critical}} \sum_{\chi_u} \prod_{i=1}^{r} \frac{1}{r^{x_i}} \xi_i^{y_i}$. This gives the coefficient of the leading term for the number of copies of $F$ in the graph formed by appending to the Turán graph $T_r(n)$ a vertex $z$ with neighborhoods of density $\boldsymbol{\xi} = (\xi_1,\ldots,\xi_r)$ with the $r$ parts. Page 5

$\pi_F := \begin{cases} \frac{\alpha_F}{|\zeta_F|}, & \text{if } \zeta_F \neq 0, \\ \infty, & \text{if } \zeta_F = 0. \end{cases}$ Page 5

$p(\rho) := \min\{P_F(\boldsymbol{\xi}) : \boldsymbol{\xi} \in \mathcal{S}_\rho\}$. Page 7

$\rho_F := \begin{cases} \inf\left\{\rho \in (0,\frac{1}{r}) : p(\rho) \leq \alpha_F \rho\right\}, & \text{if } \deg(P_F) \geq r+1, \\ \infty, & \text{if } \deg(P_F) = r. \end{cases}$ Page 7

$\hat{\rho}_F := \begin{cases} \inf\left\{\rho \in (0,\frac{1}{r}) : p(\rho) < \alpha_F \rho\right\}, & \text{if } \deg(P_F) \geq r+1, \\ \infty, & \text{if } \deg(P_F) = r. \end{cases}$ Page 7

$\mathcal{S} := \{\boldsymbol{\xi} \in \mathbb{R}^r : \forall i \in [r]\ 0 \leq \xi_i \leq 1/r\}$. Page 7

$\mathcal{S}_\rho := \{\boldsymbol{\xi} \in \mathcal{S} : \sum_{i=1}^{r} \xi_i = \rho + \frac{r-1}{r}\}$. Page 7

$t_F(n,q)$ : the smallest number of $F$-subgraphs that can be achieved by adding $q$ edges to a maximum $F$-free graph on $n$ vertices (which will be a Turán graph when $n$ is large and $F$ is color-critical). Page 1

$\mathcal{T}_r(n,q)$ : the set of graph obtained from the Turán graph $T_r(n)$ by adding $q$ edges. Page 2

$\boldsymbol{\xi} := (\xi_1,\ldots,\xi_r) \in \mathbb{R}^r$ is a vector representing the neighborhood densities of a vertex $z$ among subsets $V_1,\ldots,V_r$ of the vertex set. Page 5

$\zeta_F$ : the coefficient of $n^{f-3}$ in $\frac{\partial c}{\partial_1}(n/r,\ldots,n/r) - \frac{\partial c}{\partial_2}(n/r,\ldots,n/r)$. Roughly speaking, $\zeta_F n^{f-3}$ is the main term of the difference $c(n_1,n_2,\ldots,n_r;F) - c(n_1-1,n_2+1,\ldots,n_r;F)$. Page 5